\documentclass[11pt]{article}

\usepackage[margin=1in]{geometry}
\usepackage{amsmath,amssymb,amsthm,mathtools}
\usepackage{booktabs}
\usepackage{hyperref}
\usepackage{microtype}
\usepackage{authblk}
\newtheorem{lemma}{Lemma}
\newtheorem{theorem}{Theorem}

\theoremstyle{remark}
\newtheorem{remark}{Remark}

\newcommand{\Rop}{\operatorname{Rop}}
\newcommand{\Thi}{\operatorname{Thi}}
\newcommand{\len}{\operatorname{Len}}
\newcommand{\stick}{\operatorname{stick}}
\newcommand{\Cr}{\operatorname{Cr}}

\newcommand{\E}{E_N}
\newcommand{\length}{\operatorname{Len}}

\usepackage{graphicx}

\title{LLM-Assisted Advances in the Ropelength Problem}
\author{Alexander R. Klotz\footnote{alex.klotz@csulb.edu}}
 \affil{Department of Physics and Astronomy, California State University, Long Beach}

\begin{document}
\maketitle
\begin{abstract}
This work summarizes some updates to the ropelength problem, concerning the shortest possible contour length of a unit-thickness knot, that were derived with the assistance of ChatGPT-5.6 Luna. The first of two main results is that the universal lower bound for a knot with crossing number $C$ has been improved from $L\geq1.105C^{3/4}$ to $L\geq2.0996C^{3/4}$, achieved by further separating local and global terms in an optimization integral. The second main result is a bound relating ropelength to stick number (s), such that $L\geq\pi(s-1)$. This implies stronger lower bounds for the ropelengths of at least 41 torus knots, and a specific class of pretzel knots. Other results include an improved packing-based lower bound for $T(Q,Q)$ links and an optimized standard parameterization for torus knots. There were many problems for which the LLM was unable to render assistance, for example in reducing the numerical coefficients or scalings of known universal upper bounds.
\end{abstract}

\setcounter{tocdepth}{1}

\tableofcontents

\section{Introduction and Summary of Results}

The ropelength problem asks what the tightest possible configuration of a knot is, and expresses the ropelength $L$ as the minimum ratio of the contour length of the knot to its thickness, where thickness is defined as the radius of the smallest circle passing through any three points on the knot \cite{cantarella2002minimum}. This definition allows a normalization of the knot such that local regions do not have a radius of curvature below 1, and distant parts of the knot do not come within distance 2 of each other, allowing the knot to be surrounded by a thick tube of radius 1. Upper bounds for specific knots can be found from gradient descent, but more generally we may ask what are the upper and lower bounds for the ropelength of a knot with a given crossing number $C$. There are established universal upper and lower bounds for all knots, and these are typically far from each other and much higher and lower (respectively) than the estimated ropelengths of specific knots determined by gradient descent. There is considerable room for improvement in these bounds, and some tighter bounds for specific classes of knots (such as torus knots) have been derived.

Recently, large language models (LLMs), commonly referred to as artificial intelligence (AI) have made significant progress in producing mathematical proofs. In this manuscript, I share the results of an extended prompting session with ChatGPT-5.6 Luna. Some significant results arose out of this effort. I do not think any of them involve insights that the right humans could not have arrived at, and the most novel involves a missed connection in the literature. There are certain questions that it could not answer, and it made many mistakes in the process, requiring human oversight. The results, both positive and negative, should not be thought of as setting the end-point of the problem, but may be useful to guide further human investigators to a starting point that they can try to surpass. In the future, LLM assistance may be so ubiquitous that there is no need to point it out, but at the time of writing these norms are rapidly evolving.

The bullet points below outline the main findings in this manuscript. Several of the new bounds are plot in Fig. 1, most notably the improved Buck-Simon bound and the stick-based lower bounds for torus knots. I discuss their significance in the subsequent paragraphs. The AI-generated derivations and proofs follow. I have made comments and edits to the LLM-generated writeups where necessary. 
\begin{itemize}
    \item The numerical coefficient for the universal lower bound of the ropelength of all knots and links has been improved from 1.105 to 2.0996, such that $L\geq2.0996^{3/4}$.
    \item A lower bound for ropelength has been established by constraining it to the stick number of a knot and its ropelength, such that $L\geq\pi(s-1)$. The provides a lower bound for certain $T(P,Q), P>Q$ torus knots of $L\geq\pi(\min(4Q,2P)-1)$ that is stronger than other existing bounds for at least 41 knots. It establishes a strong linear lower bound for certain pretzel knots of $L\gtrsim (2\pi/3)C$, which is conjectured to also apply to alternating torus knots.
      \item  A somewhat tight family of configurations for $P(2,3,3,3...)$ pretzel knots was developed with $L\approx11C$, to close the sandwich opened by the stick-based lower bound.

    \item Simple geometric parameters were found for tight versions of the standard parameterization of a $T(P,Q)$ torus knot. These are comparable to the values determined from gradient descent for low crossing number,  but are not competitive as both $P$ and $Q$ increase. These are used to establish torus-specific upper bounds. 
        \item An incrementally stronger version of the packing-based lower bound for $T(Q,Q)$ links was found.
    \item An incomplete proof in a previous work from the author has been proven, establishing more firmly an upper bound for the ropelength $T(Q,Q)$ torus links. Ropelength upper bounds for $T(2Q,Q)$ and $T(3Q,Q)$ links have also been determined.
    \item A relation between Mobius energy and ropelength was extended to include links, clarifying but not proving a conjecture by the author.
  
\end{itemize}

The following query paths \textbf{did not} lead to results:
\begin{itemize}
\item No improvement to the universal upper bound of Diao-Ernst-Por-Ziegler \cite{diao2019ropelengths} was found, nor could any improved numerical parameters in Diao's lower bound be determined \cite{diao2003lower}.
\item Although it is known that the packing-based lower bound of Cantarella et al. \cite{cantarella2002minimum} is not sharp as not every component in a link can simultaneously be of minimal length, no numerical improvement could be derived that would take this into account. Similarly, the cone-piercing condition that allows lower bounds based on hexagonally packed disks could not be extended to the $T(pQ,Q)$ for $p>1$.
\item A relationship between the Milnor integrals of Brunnian links and their ropelength \cite{komendarczyk2019ropelength}  could not produce any bounds that were stronger than the standard Buck-Simon bound. Specificallly, for a family of 3-Brunnian links that includes Borromean rings and L10a140, the lower bound was found to scale as $C^{3/8}$.
\item Attempting to find a parameterization of a volume-filling torus knot was unsuccessful, despite the LLM claiming otherwise.
\end{itemize}

\begin{figure}
    \centering
    \includegraphics[width=1\linewidth]{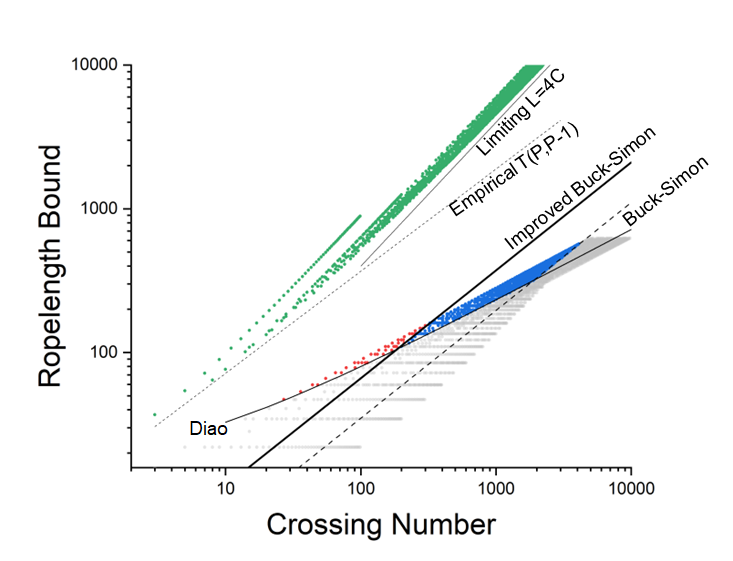}
    \caption{Ropelength bounds resulting from this work, as a function of crossing number. Thick black line shows the improved Buck-Simon bound. The lower set of points correspond to torus knots constrained from below by their stick number. The red points are a stronger lower bound than either Diao's or the improved Buck-Simon bound, while the blue points exceed the standard Buck-Simon bound (dashed line). Green points above are efficient standard torus knot parameterizations for all unique knots with $P$ and $Q$ below 100. The empirical estimates for $T(P,P-1)$ knots from Klotz and Maldonado \cite{klotz2021ropelength} and the limiting case of a standard torus of $L=4C$ are shown for comparison.}
    \label{fig:placeholder}
\end{figure}

\subsubsection*{Established bounds}

Although there are no families of knots known that scale super-linearly with crossing number, the universal upper bound is super-linear with $L\leq a_u C\log^5C$, with $a_u$ unknown except that it must be greater than 78.5 to include the Hopf link \cite{diao2019ropelengths}. There are stronger upper \cite{hong2014minimum} and lower \cite{Schoen2026EndpointTangents, diao2003lower} bounds at finite crossing numbers. The current upper bounds are

\[
\operatorname{L}(C)
\leq
\begin{cases}
\displaystyle
\left(1+\frac{2}{\pi}\right)C^2
+\left(\frac{8}{\pi}+2+\pi\right)C
+\frac{8}{\pi}+\frac{4\pi}{3},
& 3\leq C\leq 10,\\[2ex]
\displaystyle
\frac{3}{2}C^2+(\pi+6)C+2\pi+\frac{1}{2},
& 11\leq C\leq 23,\\[2ex]
\displaystyle
\frac{5}{4}C^2+
\left(\frac{\pi}{2}+13\right)C
+\pi+\frac{55}{4},
& 24\leq C\leq 47595,\\[2ex]
\displaystyle
272C^{3/2}+168C+44\sqrt{C}+22,
& 47596\leq C\leq 10^{14},\\[2ex]
\displaystyle
78.5\,C(\ln C)^5,
& C\geq 10^{14}.
\end{cases}
\]

Since every prime knot (but not link) with up to 12 crossings has been tightened using gradient descent \cite{ashton2011knot,klotz2024ropelength}, an empirical upper bound in that range may be written as $L<9C+8$. The exact crossover to the final upper bound depends on certain unknown leading terms and prefactors. The current best lower bounds are

\[
\operatorname{L}(C)
\geq
\begin{cases}
31.51217852,
& 3\leq C\leq 8 ,\\[2ex]
\displaystyle
\frac{1}{2}
\left(
17.334+
\sqrt{17.334^2+64\pi C}
\right),
& 9\leq C\leq 1896,\\[3ex]
\displaystyle
1.105\,C^{3/4},
& C\geq 1897.
\end{cases}
\]

The first is a bound on non-trivial knots in general \cite{Schoen2026EndpointTangents}, the second is known as the Diao lower bound \cite{diao2003lower} and the third as the Buck-Simon bound \cite{buck1999thickness}. It has been estimated from aforementioned numerical studies that all 10-crossing knots have ropelengths between 71 and 93 \cite{ashton2011knot}, but the universal lower and upper bounds for $C=10$ are 32.7 and 247.3. The factor of 7 gap between upper and lower compared to the factor of 1.3 gap in the numerical results is a source of motivation to derive stronger bounds.


\subsubsection*{An improved universal lower bound}

The most significant result is the improved numerical constant in the Buck-Simon universal lower bound \cite{buck1999thickness}. The Buck-Simon bound is based on a hierarchy of energy integral functional-based inequalities, ultimately constraining the ropelength to a ``normal'' (as in orthogonal) energy, to a ``symmetric'' energy, to the average crossing number, to the minimal crossing number. A key step is maximizing the normal energy, which is constrained by the unattainable case arising when the knot is maximally packed within a spherical annulus surrounding a given point. The new solution improves the packing argument by further separating the local and global contributions to the energy integral and reducing a numerical constaint associated with the packing from 10.67 (rounded to 11 by Buck and Simon) to 4.67. The ropelength lower bound is

\begin{equation}
    L\geq\left(\frac{4\pi}{4.674}\right)^{3/4}C^{3/4}\approx2.0996C^{3/4}>\left(\frac{4\pi}{11}\right)^{3/4}C^{3/4}\approx1.105C^{3/4}.
\end{equation}

This is almost a factor-of-two improvement in the numerical coefficient. Diao's lower bound \cite{diao2003lower}, which is stronger than the Buck-Simon bound for knots with up to 1896 crossings, is exceeded by this improved version at 187 crossings. However, it is still likely far below the ``true'' lower bound for all knots and links, as the smallest values of $L/C^{3/4}$ that have been observed from gradient descent are above 10.

\subsubsection*{A connection between stick number and ropelength}

The most novel finding to come out of this, in the author's opinion, is a neglected relation between stick number and ropelength. This arose from the prompt ``Do a review of more recent literature on knot energies, besides Mobius and the Buck cross-energy, and see if there are constraints on ropelength that have not been appreciated or stated explicitly.''  It follows an argument from Litherland et al. \cite{litherland1999thickness} which constrains the number of straight edges that a thick tubular knot can be broken into, visualized in Fig. 2. The novel contribution of ChatGPT was to relate this minimum number to the stick number ($s$) of knots, which sets a lower bound for the number of line segments in a polygonal knot \cite{CantarellaRechnitzerSchumacherShonkwiler2026} and is constrained by other analytical techniques and numerical methods. If a knot has  unit thickness then the relation is

\begin{equation}
    s\leq\frac{L}{\pi}+1\longrightarrow L\geq\pi(s-1.
\end{equation}

The universal lower bound for stick number depends on the square root of crossing number \cite{CantarellaRechnitzerSchumacherShonkwiler2026}, so this ropelength bound will not be globally competitive. However, there are cases in which it can be. The stick number of torus knots (with $P>Q$) is  $s(T(P,Q))=\min(2P,4Q)$ \cite{Jin1997}, which implies a ropelength lower bound of

\begin{equation}
    L(T(P,Q)\geq\pi(\min(4P,2Q)-1).
\end{equation}

This is stronger than all known bounds for 41 torus knots with between 27 and 300 crossings. Using the standard Buck-Simon bound, the torus stick bound is superior for 428 knots. These knots are listed in the section 3.5.

There is a conjectured stick number for alternating torus knots ($T(P,2)$), which was observed numerically to hold for $P$ from 3 to 31, of the form $s(P)=\lfloor\frac{2}{3}(P+1)\rfloor+4$. If this were a true lower bound, it would imply a ropelength for these knots of

\begin{equation}
L(T(p,2))\geq\pi\lfloor \frac{2(p+1)}{3}\rfloor+3\approx \frac{2}{3}\pi p
\end{equation}
This ropelength lower bound is asymptotically linear with a slope of $2.09$, considerably stronger than the current bound of 1/56 \cite{diao2022ropelength}, but still much less than the best known upper bound of 7.31 \cite{kim2024efficiency}.

\begin{figure}
    \centering
    \includegraphics[width=1\linewidth]{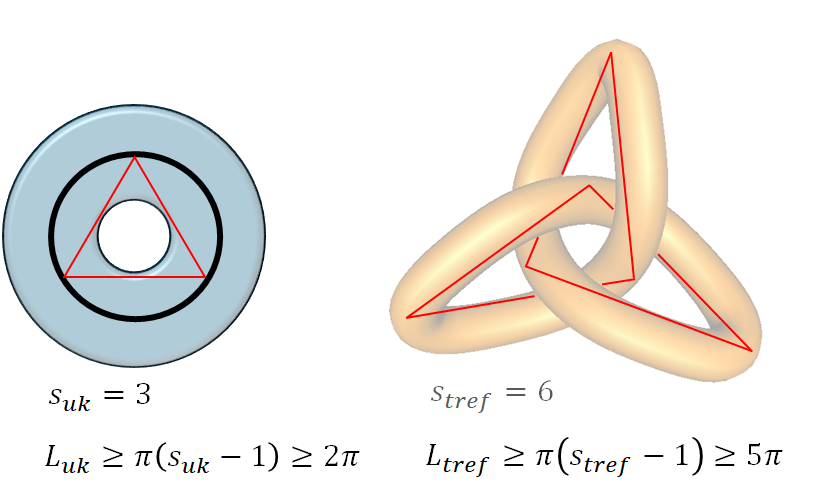}
    \caption{Schematic showing a smooth tube surrounding an unknot (blue) or trefoil (orange) and how the central curve may be replaced with piecewise linear segments that do not pierce the tube. This may be used to constrain the stick number to the ropelength. This image was human-generated.}
    \label{fig:stick}
\end{figure}

The stick relation is able to provide a strong lower bound for the family of $P(2,3,3,3...)$ pretzel knots described by Cantarella et al. \cite{CantarellaRechnitzerSchumacherShonkwiler2026} which have a known lower bound for their bridge and stick number. Together these imply a ropelength of

\begin{equation}
    L_{pretz}\geq \frac{2\pi}{3}\,C+\frac{5\pi}{3}.
\end{equation}

This is a factor of 75 greater than the lower bound of $C/36$ established by Diao and Ernst \cite{diao2025}. Other families of pretzel knots as well as certain non-pretzel Montesinos knots can have a strong linear slope from a similar argument, but the $P(2,3,3,3...)$ family is the strongest. Interestingly, this could have been realized from the bridge number alone, which already sets a bound on ropelength \cite{DiaoErnstThistlethwaite2003}, but the connection was not made in the literature.

\subsubsection*{Tight configurations of pretzel knots}
Given the strong lower bounds for the $P(2,3,3,3...)$ pretzel knot family, I directed the LLM to construct tight configurations for a strong upper bound as well. It was unable to do this autonomously, creating Cartesian instances of a braid representation that were hihgly inefficient. I then directed it to construct them out of a circular arrangement of efficient double helices, connected at the top and bottom by circular arcs. It generally followed my instructions and found a way to connect adjacent helices with Bezier curves. Based on my visual inspection, I suggested improvements to reduce the overall lengths, which eventually lead to configurations with a differential ropelength of approximately 11 per crossing. Inspecting such a configuration with 32 crossings (Figure 3), there is clear slack that can be removed and I consider this a good first step but not a complete success. The 32-crossing instance of a pretzel knot had an initial length of 363  which was somewhat reduced to 319 by Ridgerunner, a factor of 4.4 above the lower bound of 72.3. This minimization likely had not converged. Combined with the stick-based lower bound, these configurations establish knots with an upper and lower ropelength bound that are within a factor of 6 of each other.

\begin{figure}
    \centering
    \includegraphics[width=0.5\linewidth]{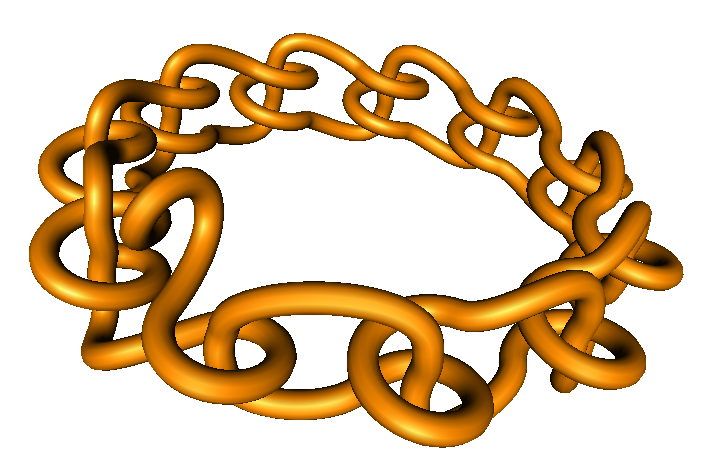}
    \caption{Instance of a tight $P(2,3,3,3,3,3,3,3,3,3,3)$ pretzel knot.}
    \label{fig:placeholder}
\end{figure}

\subsubsection*{Clarification and mild improvement of $T(Q,Q)$ and $T(2Q,Q)$ bounds}

In a recent paper \cite{klotz2026tight}, I derived strong upper and lower bounds for the ropelength of $T(Q,Q)$ torus links. The upper bounds were based on constructing two densely packed tori of concentric helices and linking them through each other. My derivation for the upper bound relied on an approximation to the no-overlap constraint inequality for $N$ helices wrapped around a cylinder. My proof that it was non-overlapping was incomplete because I could only show that constraints weren't violated in a third-order expansion of the inequality. The LLM has produced a proof of my conjecture that does not rely on a cubic approximation. The upper bounds also relied on a small numerical correction (on the order 1.001) for the arc length of a toroidal helix, for which the LLM has found an approximation.

While I worked out the upper bounds for $T(Q,Q)$ knots as $L\leq11.68C^{3/4}$, I mentioned that a double linking of two $T(2Q,Q)$ or $T(3Q,Q)$ links would have stronger upper bounds, although they would no longer be torus links. I did not work out the bounds explicitly. The LLM determined that they would have coefficients of 11.21 and 11.77 respectively, the $2Q$ version having a stronger bound (but also a weaker lower bound).

The lower bounds were based on applying the Wegner disk packing formula \cite{wegner1986uber} to the Cantarella bound relating ropelength to the convex hulls around disks \cite{cantarella2002minimum}. Specifically, the bound states that the convex hull around $N$ unit disks cannot be shorter than the perimeter of a circle that has the same area as the minimum-area hexagon-shaped hexagonal packing of $N$-disks, known as a Wegner packing. The LLM was able to improve upon this slightly by replacing the circular isoperimetric bound with a regular $N$-gon isoperimetric bound. The updated lower bound is  

\begin{equation}
 L_{T(Q,Q)}
 \geq
 Q\left[
 2\pi+\sqrt{
 4\pi A_W(Q-1)+\frac{4\pi^2}{3}}
 \right].
 \label{eq:improved-tqq-boun2d}
\end{equation}

where $A_W(Q-1)$ is the area of a Wegner packing of $Q-1$ disks. This has the effect of increasing the lower bound of the $T(10,10)$ link from 818 to 828, and the effect gets smaller as $1/Q$. It is still unproven that the hexagonal packing bounds can apply to $T(pQ,Q)$ links for $p>1$ due to a cone-piercing condition established by Cantarella et al. \cite{cantarella2002minimum}. The LLM was unable to prove that this cone-piercing condition is satisfied.

\subsubsection*{Tight standard torus knot parameterization}

This path was motivated by a desire to be able to write down a simple equation for a tight knot that does not rely on opaque numerical parameters. The LLM was prompted to find ``nice'' geometric principles for generating tight torus knots based on the standard harmonic parameterization, where nice refers to simple functions of $P$ and $Q$.  I directed the LLM to optimize this as a function of one to two free parameters: the ratio of the minor to major torus radii ($\rho$), and a scaling of the z-coordinates ($\lambda$) such that

\begin{equation}
 \gamma_{P,Q}(t)
 =
 \bigl((\rho+\cos Pt)\cos Qt,\,
       (\rho+\cos Pt)\sin Qt,\,
       \lambda\sin Pt\bigr),
 \qquad 0\leq t\leq 2\pi .
\label{eq:standard-torus1}
\end{equation}

It arrived at the following optimization (with $P>Q$) for $\rho$:

\begin{equation}
    \rho=1+\frac{P}{Q}
\end{equation}

and two reasonable solutions for $\lambda$:

\begin{equation}
    \lambda_1=\sqrt{1+\frac{Q^2}{P^2}} \ \ \ \ \ \ \ \ \ \ \ \lambda_2=\csc\left(\frac{\pi Q}{P}\right)
\end{equation}

Typically, the square root optimization is useful when $P>7Q$ and the choice is non-trivial below that. $\lambda=1$ works best for $Q=2$. However, the optimal z-scaling will reduce ropelength from radial-only by at most 9\%. Figure 4 shows the improvement to radial-only scaling as a function of $P/Q$, and which optimization dominates.  

The green points in Figure 1 show the ropelengths produced by this parameterization. These configurations are known not to be competitive, we showed previously \cite{thompson2025ropelength} that $Q$ helices wrapped around a single cylinder have an asymptotic linear ropelength of $4C$, and before that showed that tight torus knots have a ropelength that grows empirically as $L\approx14C^{0.71}$ \cite{klotz2021ropelength}. Although in a previous work we found that the ropelength-per-crossing of torus knots when $P>>Q$ can be constructed near $7.6/\sqrt{Q}$, we did not have a generation solution for $P\gtrsim Q$.

The LLM was able to derive an exact expression for the ropelength of these construtions which relies on a hypergeometric integral and the solution of a transcendental equation, which is not very useful. However, this can be approximated to find the linear slope of the ropelength of a $T(P,Q)$ knot for fixed $P$ as $Q$ varies. This yields an approximate ropelength slope for fixed $Q$ of

\begin{equation}
    M_Q\approx4Q+\frac{\pi^2}{12Q}
\end{equation}

\begin{figure}
    \centering
    \includegraphics[width=1\linewidth]{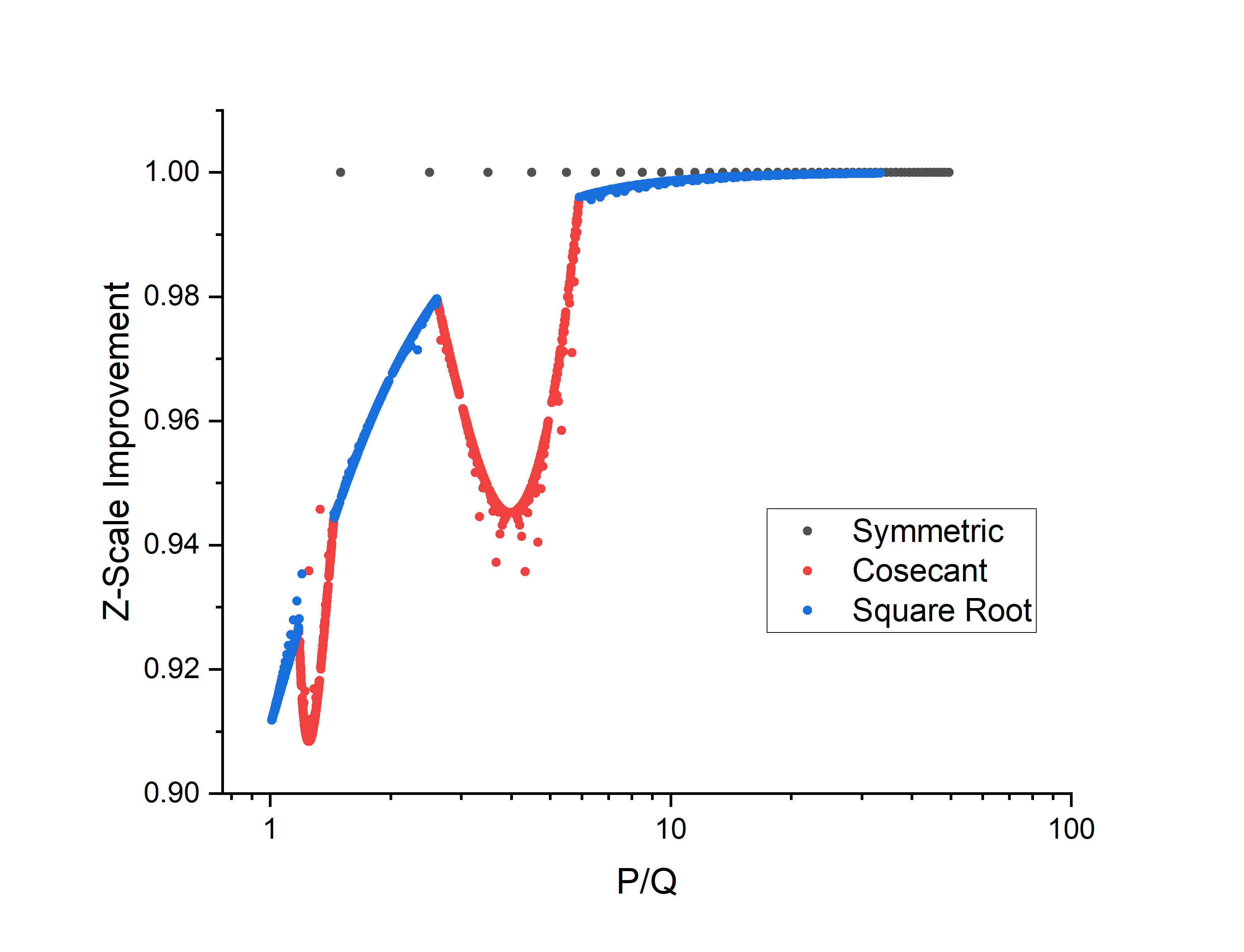}
    \caption{The ratio between the ropelengths of two-parameter optimized standard parameterization torus knots, to those with only radial optimization, for all knots with $Q<P\leq100$. The colors indicate which Z-scale factor is optimal. Radial-only optimization is optimal for $T(P,2)$ knots.}
    \label{fig:placeholder}
\end{figure}

Because this approximation is an underestimate also found a linear upper bound for the ropelength of torus knots. Although this derivation explicitly included the curvature and self-intersection criteria, it refuses to state that this serves as an upper bound.

\begin{equation}
  L_u(P,Q)
 =
 2\pi QP+\pi Q^2+O(P^{-1}).
\label{eq:UPQ-two-denominator2}
\end{equation}

None of these are competitive with my previous work on concentric helices \cite{thompson2025ropelength}, but these represent solutions for closed single-component knots, which my previous work did not realize. These constructions can be thought of as analogous for the general $Q$ to the Olsen-Bohr solution for $T(P,2)$ knots \cite{olsen2012geometry}, which was the first of several gradually improved optimizations.

\subsubsection*{Clarification on Mobius-based bounds for links}
On a minor note, there exists a bound from Rawdon and Simon \cite{rawdon2002mobius} between the minimal Mobius energy of knots and their ropelength. Because we do not have good lower bounds for Mobius energy, this has not been useful for bounding ropelength, but I have previously used this bound in tandem with some other (unproven) assumptions to make a conjectured improvement to the universal lower bound \cite{klotz2025geometric}. This conjecture relies on the Rawdon-Simon bound applying to links. The LLM showed that

\begin{equation}
    E_M\leq4.64L^{4/3}
\end{equation}

replacing the 4.56 used for knots. This also implies that my conjectured bound of $3.22C^{3/4}$ would actually have the coefficient 3.19. However, my conjecture is unproven. It also relies on an unresolved choice of convention for the Mobius energy of a link, if a different convention is used than the coefficient is reduced to 2.02, which is no longer greater than the improved Buck-Simon bound.

\subsection*{Attempted LLM-based ``peer'' review}

After compiling an initial version of this manuscript, I uploaded the PDF to Microsoft Copilot and asked for a peer referee-style report. It had some useful comments, then lied to me twice about compiling them into a PDF, instead creating much lighter summaries of the report. Most of the criticisms focused on the fact that numerical parameters were not certified. I showed this report to ChatGPT which claimed it would make these improvements (it could not, as it was not actually writing the paper). In response, I had ChatGPT generate some supporting numerical evidence for its claims and a Python script to replicate the numerical integral values.

\subsection*{Final Human Comments}

Note that most of the following text has been generated by the LLM. Each section may have a unique style that differs from other sections of the manuscript, and basic definitions may be repeated. I have added a few comments in bold throughout and made some minor edits. In the future, it may be standard to have papers that are co-written or fully written by LLMs and it will be redundant to state that, but in 2026 the norms have not converged. It should be noted that ChatGPT and Copilot are not faithful reporters of their own outputs. Requested LaTeX files are often different than the chat outputs they are requested to reproduce, which can be a source of errors. I recommend to anyone attempting similar LLM-assisted research to have a LaTeX text of each response generated alongside the response.

\bibliographystyle{unsrt}
\bibliography{knotrefs}

\section{A stronger universal Buck-Simon lower bound}

\textbf{Author's note: the original Buck-Simon derivation uses the following chain of inequalities:}
\begin{equation}
    11L^{4/3}\geq E_N\geq E_s\geq4\pi\langle C\rangle\geq4\pi C
\end{equation}
\textbf{where $E_N$ is the normal knot energy, $E_S$ is the symmetric knot energy, and $\langle C\rangle$ is the average crossing number. These are all functionals in the form of double integrals over the knot and are defined in the Buck-Simon paper \cite{buck1999thickness}. Importantly, the factor of 11 comes from deriving a bound on the normal energy that assumes the knot is maximally packed into spherical shells. The argument below attempts to improve the numerical parameters from the packing argument.}

\subsection{Normalization and notation}

Let $K\subset\mathbb{R}^3$ be a smooth knot with thickness (tube radius)
$\Delta>0$.  Its ropelength is
\[
    \Rop(K)=\frac{\length(K)}{\Delta}.
\]
Because both length and thickness scale linearly, we may rescale so that
$\Delta=1$.  In this normalization write
\[
    L=\length(K)=\Rop(K).
\]
Thus $K$ has a non-overlapping unit-radius normal tube.

For $x\in K$, let $T_x$ denote the unit tangent vector at $x$.  If
$x=\gamma(s)$ for an arclength parametrization $\gamma$, then
\[
    T_x=\gamma'(s).
\]

The normal energy is
\begin{equation}
\label{eq:normal-energy}
    \E(K)=
    \int_K\int_K
    \frac{\bigl|T_x\times(y-x)\bigr|^2}{|y-x|^4}\,ds_y\,ds_x.
\end{equation}
Since $|T_x|=1$,
\[
    |T_x\times(y-x)|^2
    =|y-x|^2-\bigl(T_x\cdot(y-x)\bigr)^2.
\]
In particular,
\begin{equation}
\label{eq:integrand-crude}
    \frac{|T_x\times(y-x)|^2}{|y-x|^4}
    \le \frac{1}{|y-x|^2}.
\end{equation}

\subsection{The topological input}

The Buck--Simon energy--crossing estimate is
\begin{equation}
\label{eq:BS-lower}
    \E(K)\ge 4\pi\,c(K).
\end{equation}
Consequently, an upper estimate of the form
\begin{equation}
\label{eq:energy-upper-form}
    \E(K)\le A L^{4/3}
\end{equation}
implies
\begin{equation}
\label{eq:ropelength-from-energy}
    L\ge \left(\frac{4\pi}{A}\right)^{3/4}c(K)^{3/4}.
\end{equation}

The task is therefore to obtain the smallest possible universal coefficient
$A$ from a geometric upper estimate for $\E(K)$.

\subsection{A fixed-point decomposition}

Fix $x\in K$.  Let $K_{\rm loc}(x)$ consist of the two arclength pieces
starting at $x$ and extending distance $3$ in the two directions along $K$.
The remaining part is denoted by
\[
    K_{\rm far}(x)=K\setminus K_{\rm loc}(x).
\]
For the knots under consideration, $L\ge 6$, and the local piece has length
$6$.

Define
\[
    I(x)=\int_K
    \frac{|T_x\times(y-x)|^2}{|y-x|^4}\,ds_y.
\]
Then
\[
    \E(K)=\int_K I(x)\,ds_x.
\]
We write
\[
    I(x)=I_{\rm loc}(x)+I_{\rm far}(x).
\]

\subsection{The local contribution}

The standard Buck--Simon circle comparison gives the pointwise estimate
\begin{equation}
\label{eq:local-integrand}
    \frac{|T_x\times(y-x)|^2}{|y-x|^4}\le \frac14
\end{equation}
for points in the two local arclength pieces.  Since their total length is
$6$, equation \eqref{eq:local-integrand} gives
\begin{equation}
\label{eq:local-bound}
    I_{\rm loc}(x)\le \frac{6}{4}=\frac32.
\end{equation}

\subsection{The far contribution}

The Schur comparison used in the Buck--Simon argument gives the separation
scale
\begin{equation}
\label{eq:d-definition}
    |x-y|\ge d:=2\sin\left(\frac32\right)
\end{equation}
for the far points.  Numerically,
\[
    d=1.994989973208\ldots .
\]
Because the tube radius is one, the unit tube about the far portion stays
outside the ball of radius
\begin{equation}
\label{eq:b-definition}
    b=d-1=0.994989973208\ldots
\end{equation}
centered at $x$.

For $R\ge d$, define the cumulative far length
\[
    M_x(R)=\length\{y\in K_{\rm far}(x):|x-y|\le R\}.
\]
The unit tube about this part of the curve is contained in the spherical
annulus
\[
    b\le |z-x|\le R+1.
\]
The annulus has volume
\[
    \frac{4\pi}{3}\bigl((R+1)^3-b^3\bigr).
\]
Dividing by the cross-sectional area $\pi$ gives the basic cumulative
estimate
\begin{equation}
\label{eq:F-definition}
    M_x(R)\le F(R),
    \qquad
    F(R)=\frac43\bigl((R+1)^3-b^3\bigr).
\end{equation}

\subsection{The improvement from the local tube}

The basic estimate \eqref{eq:F-definition} treats the annulus as if it were
available entirely to the far tube.  However, the local tube also occupies
part of that annulus.

The unit tube about the local arc has volume $6\pi$.  At most the volume
$4\pi b^3/3$ of the ball $B(x,b)$ can lie inside that ball.  Therefore at
least
\[
    6\pi-\frac{4\pi}{3}b^3=\pi\delta
\]
of the local tube lies in the annulus $b\le |z-x|\le R+1$, provided $R\ge3$.
Here
\begin{equation}
\label{eq:delta-definition}
    \delta=6-\frac43b^3
    =4.686606540032\ldots .
\end{equation}

The local and far unit tubes are disjoint.  Consequently, for $R\ge3$,
\begin{equation}
\label{eq:improved-F}
    M_x(R)\le F(R)-\delta.
\end{equation}

For $d\le R\le3$, monotonicity gives
\[
    M_x(R)\le M_x(3)\le F(3)-\delta
    =\frac{238}{3}.
\]
Combining these estimates, define
\begin{equation}
\label{eq:H-definition}
H(R)=
\begin{cases}
0, & R<d,\\[1mm]
\min\left\{F(R),\dfrac{238}{3}\right\}, & d\le R\le3,\\[3mm]
F(R)-\delta, & R\ge3.
\end{cases}
\end{equation}
Then, if
\[
    S=\length(K_{\rm far}(x))=L-6,
\]
we have the effective cumulative bound
\begin{equation}
\label{eq:cumulative-bound}
    M_x(R)\le \min\{S,H(R)\}.
\end{equation}

\subsection{Layer-cake estimate for the far energy}

Using \eqref{eq:integrand-crude},
\[
    I_{\rm far}(x)
    \le \int_{K_{\rm far}(x)}\frac{ds_y}{|x-y|^2}.
\]
For $r>0$,
\[
    \frac{1}{r^2}=\int_r^\infty \frac{2}{R^3}\,dR.
\]
Applying this identity and interchanging the order of integration gives
\begin{align}
    I_{\rm far}(x)
    &\le \int_d^\infty \frac{2M_x(R)}{R^3}\,dR \\
    &\le G(S),
\end{align}
where
\begin{equation}
\label{eq:G-definition}
    G(S)=\int_d^\infty
    \frac{2\min\{S,H(R)\}}{R^3}\,dR.
\end{equation}

Combining \eqref{eq:local-bound} and \eqref{eq:G-definition},
\[
    I(x)\le \frac32+G(L-6).
\]
This bound is independent of $x$, so integration over $K$ yields
\begin{equation}
\label{eq:energy-upper}
    \E(K)\le L\left(\frac32+G(L-6)\right).
\end{equation}

Define the dimensionless objective
\begin{equation}
\label{eq:Q-definition}
    Q(L)=\frac{\frac32+G(L-6)}{L^{1/3}}.
\end{equation}
Then \eqref{eq:energy-upper} becomes
\begin{equation}
\label{eq:energy-Q}
    \E(K)\le Q(L)L^{4/3}.
\end{equation}

\subsection{Numerical maximization}

The function $F$ reaches the plateau value $238/3$ at
\[
    R_0=2.925388714960\ldots .
\]
On the branch containing the maximum, the active cutoff radius $R$ is
related to $L$ by
\begin{equation}
\label{eq:L-R-relation}
    L=F(R)-\delta+6,
    \qquad R\ge3.
\end{equation}

Numerical quadrature of \eqref{eq:G-definition}, followed by one-dimensional
maximization of \eqref{eq:Q-definition}, gives
\begin{equation}
\label{eq:numerical-max}
    R_*=7.463029734\ldots,
    \qquad
    L_*=808.1953302\ldots,
\end{equation}
and
\begin{equation}
\label{eq:Qmax}
    Q(L_*)=4.673424743\ldots .
\end{equation}
Thus the conservative rounded estimate
\begin{equation}
\label{eq:Q-conservative}
    Q(L)<4.674
\end{equation}
gives
\begin{equation}
\label{eq:energy-conservative}
    \E(K)\le4.674\,L^{4/3}.
\end{equation}

\textbf{Author's note: to shed light on the asserted numerical maximum, I directed the LLM to produce a plot of this function showing the maximum (Fig. 5). The plot itself is based on a purported numerical integration which is hard to validate, and I had the LLM generate a Python script to compute this, which is attached in the ancillary material.}

\begin{figure}
    \centering
    \includegraphics[width=1\linewidth]{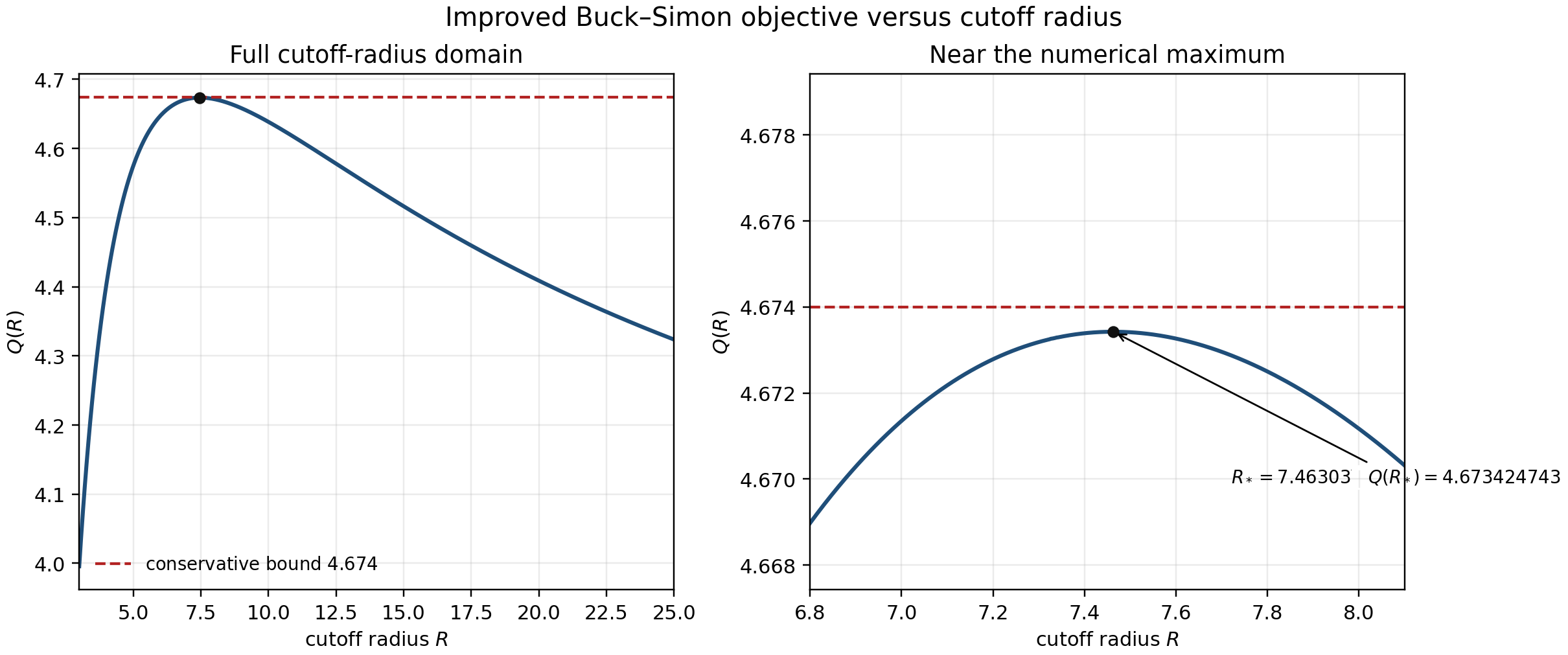}
    \caption{LLM-generated plot showing maximum of energy function. A small graphical glitch on the right figure was my attempt at correcting a typesetting issue in the original version.}
    \label{fig:energymax}
\end{figure}

\subsection{Conversion to a ropelength lower bound}

Combining \eqref{eq:BS-lower} and \eqref{eq:energy-conservative},
\[
    4\pi c(K)\le4.674\,L^{4/3}.
\]
Solving for $L$ gives
\begin{equation}
\label{eq:final-bound}
    \boxed{
    \Rop(K)=L
    \ge
    \left(\frac{4\pi}{4.674}\right)^{3/4}c(K)^{3/4}
    =2.0996208\ldots\,c(K)^{3/4}.
    }
\end{equation}
Rounding downward gives the stated form
\[
    \boxed{\Rop(K)\ge2.0996\,c(K)^{3/4}.}
\]

If the unrounded numerical maximum in \eqref{eq:Qmax} were certified as a
global maximum, it would give the slightly larger numerical coefficient
\[
    \left(\frac{4\pi}{4.673424743\ldots}\right)^{3/4}
    =2.0998146\ldots .
\]
The value $2.0996$ is preferable for a conservative statement.

\subsection{Status of the argument}

There are three logically distinct ingredients:
\begin{enumerate}
    \item the Buck--Simon topological inequality \eqref{eq:BS-lower};
    \item the local and far-field tube estimates, including the separation
          input \eqref{eq:d-definition};
    \item the numerical maximization of the explicit one-variable function
          \eqref{eq:Q-definition}.
\end{enumerate}

The improvement over the older $4.724$ energy coefficient comes from the
subtraction of $\delta$ in \eqref{eq:improved-F}: the local tube already
occupies part of the annular volume that would otherwise be assigned to the
far tube.  The numerical value $4.673424743\ldots$ is obtained by ordinary
high-precision quadrature and optimization.  A completely formal published
proof of the numerical maximum would additionally certify the quadrature and
the sign of the derivative on the entire admissible interval, for example by
interval arithmetic.

The geometric separation estimate used in the far-field step must likewise
be stated with exactly the hypotheses under which it is valid.  Thus the
derivation above should be read as a self-contained proof conditional on the
standard Buck--Simon separation lemma and on the unit-tube packing
properties used in their argument.  The volume subtraction itself is
elementary once those inputs are available.

\textbf{Author's note: the 4.724 mentioned above refers to an earlier finding by ChatGPT leading to a coefficient of 2.0829, and was the source of many queries, which is why the current result is being compared to it.}

\subsection{A conditional improvement from hexagonal tube packing}

A substantially larger improvement would follow from a universal
finite-shell packing estimate.  Suppose that there exists a constant
\(\eta<1\) such that every relevant unit-tube packing satisfies
\[
M_x(r)
 \le
\eta\,F_0(r),
\qquad
F_0(r)
 :=
\frac43\left((r+1)^3-1\right).
\tag{13}
\]
The conjectural hexagonal value is
\[
\eta=\frac{\pi}{\sqrt{12}}
 =0.906899682\ldots .
\tag{14}
\]

Using the original local split at arc distance \(\pi\), define
\[
G_\eta(S)
 :=
 \int_2^\infty
 \frac{2\min\{S,\eta F_0(r)\}}{r^3}\,dr.
\tag{15}
\]
The corresponding energy estimate is
\[
E_N(K)
 \le
 L\left(\frac{\pi}{2}+G_\eta(L-2\pi)\right).
\tag{16}
\]

Numerically,
\[
\sup_{L\ge2\pi}
\frac{\frac{\pi}{2}+G_\eta(L-2\pi)}{L^{1/3}}
 =
4.44131\ldots .
\]
Consequently, assumption (13) with the value (14) would yield
\[
\operatorname{Rop}(K)
\ge
\left(\frac{4\pi}{4.44131}\right)^{3/4}
\operatorname{Cr}(K)^{3/4},
\]
or
\[
\boxed{
\operatorname{Rop}(K)
\ge
2.1816\,\operatorname{Cr}(K)^{3/4}.
}
\tag{17}
\]

This \(2.18\) coefficient is conditional.  The required finite-shell
density estimate is stronger than the currently established general
tube-packing results.  In particular, the conjectural status of the
hexagonal-density bound for arbitrary thick links must be resolved
before (17) becomes an unconditional ropelength theorem.

\section{Constraining ropelength to stick number}

\subsection{Definitions and normalization}

Let $\gamma:S^1\to\mathbb R^3$ be an embedded closed curve.  Its thickness is
denoted by $\Thi(\gamma)$.  For a $C^{1,1}$ curve this can be characterized as
the largest radius of an embedded normal tube; equivalently, it is controlled
by both the radius of curvature and the relevant doubly-critical self-distance.
The ropelength of a particular representative is
\[
    \mathcal R(\gamma)
    := \frac{\len(\gamma)}{\Thi(\gamma)},
\]
and the ropelength of a knot type $K$ is
\[
    \Rop(K):=\inf_{\gamma\in K}\mathcal R(\gamma).
\]

The formulas below use the unit-radius convention.  Thus, after rescaling any
representative, we may assume $\Thi(\gamma)=1$ and then
\[
    \mathcal R(\gamma)=\len(\gamma).
\]

\subsection{The polygonal approximation input}

The following is the edge-number estimate proved in the thickness theory of
Litherland, Simon, Durumeric, and Rawdon \cite{litherland1999thickness}.  It is the precise bridge between
ropelength and stick number used here.

\begin{lemma}[Thick-curve polygonal approximation]
Let $\gamma$ be an embedded closed curve with length $\ell$ and thickness
$\tau>0$.  There is a polygonal knot $P$ ambient isotopic to $\gamma$ with
\[
    N(P)\leq \frac{\ell}{\pi\tau}+1
\]
edges.
\end{lemma}

\begin{proof}[Geometric mechanism]
Partition $\gamma$ into consecutive subarcs whose arclengths are at most
$\pi\tau$, and join consecutive partition points by straight chords.  The
thickness condition gives the local tube and restricted-bending estimates
needed to ensure that each such subarc can be replaced by its chord without
creating a self-intersection or changing the ambient isotopy class.  The
resulting inscribed polygon is therefore isotopic to $\gamma$.

If $N$ is the number of subarcs, then the partition can be chosen with
\[
    N\leq \frac{\ell}{\pi\tau}+1.
\]
The polygon has $N$ edges, which proves the claim.
\end{proof}

\begin{remark}
The nontrivial part of the lemma is the isotopy assertion.  It is not a
consequence of arclength partitioning alone: an arbitrary curve can have a
short subarc whose chord is topologically unsafe.  The thickness hypothesis is
what supplies the excluded-volume and bending control.  Accordingly, the
lemma should be cited as the polygonal-approximation theorem rather than
presented as an elementary consequence of the definition of thickness.
\end{remark}

\subsection{Derivation of the stick-based ropelength bound}

Let $s(K)=\stick(K)$ denote the minimum number of edges in any polygonal
representative of $K$.  Applying the lemma to an arbitrary representative
$\gamma\in K$ gives
\[
    s(K)\leq \frac{\len(\gamma)}{\pi\,\Thi(\gamma)}+1.
\]
Rearranging,
\[
    \frac{\len(\gamma)}{\Thi(\gamma)}
    \geq \pi\bigl(s(K)-1\bigr).
\]
Since this holds for every representative of $K$, taking the infimum yields
the main result.

\begin{theorem}[Stick-based ropelength lower bound]
For every nontrivial knot type $K$,
\[
    \boxed{\Rop(K)\geq \pi\bigl(\stick(K)-1\bigr).}
\]
More generally, if a rigorous lower bound $s_0\leq \stick(K)$ is known, then
\[
    \Rop(K)\geq \pi(s_0-1).
\]
\end{theorem}

\begin{proof}
The inequality in the preceding display holds for every representative
$\gamma\in K$.  Taking the infimum over $\gamma$ proves the result.
\end{proof}

\begin{remark}[Exactness requirement]
To use the theorem numerically, one needs an exact stick number or a rigorous
lower bound for it.  A polygonal construction with $s$ edges proves only
$\stick(K)\leq s$ and therefore does not, by itself, imply
$\Rop(K)\geq\pi(s-1)$.  The latter inference is valid when the construction is
accompanied by a proof that no polygon with fewer edges represents $K$.
\end{remark}

\subsection{Torus-knot lower bounds}

Write $T(q,p)$ with $2\leq p<q$ and $\gcd(p,q)=1$.  Its crossing number is
\[
    \Cr\bigl(T(q,p)\bigr)=q(p-1).
\]

In the ranges where the torus-knot stick-number theorem \cite{Jin1997} gives
\[
\stick\bigl(T(q,p)\bigr)=
\begin{cases}
2q, & p<q<2p,\\[3pt]
4p, & 2p<q<3p,
\end{cases}
\]
the theorem above gives
\[
\Rop\bigl(T(q,p)\bigr)\geq
\begin{cases}
\pi(2q-1), & p<q<2p,\\[3pt]
\pi(4p-1), & 2p<q<3p.
\end{cases}
\]


\subsection{Torus knots for which the stick bound is strongest}

Write $T(q,p)$ with $2\leq p<q$ and $\gcd(p,q)=1$.  Its crossing
number is
\[
    C(T(q,p))=q(p-1).
\]
In the ranges covered by the exact stick-number theorem,
\[
\operatorname{stick}(T(q,p))=
\begin{cases}
2q,&p<q<2p,\\
4p,&2p<q<3p.
\end{cases}
\]
Consequently, the stick-derived ropelength estimate is
\[
R_{\mathrm{stick}}(T(q,p))
   =\pi\bigl(\operatorname{stick}(T(q,p))-1\bigr).
\]

For comparison, define Diao's bound by
\[
D(C):=\frac{17.334+\sqrt{17.334^2+64\pi C}}{2},
\]
and define
\[
B_{2.0996}(C):=2.0996\,C^{3/4},
\qquad
B_{1.105}(C):=1.105\,C^{3/4}.
\]
All lists below impose $C(T(q,p))<1850$.

\subsubsection*{Knots beating both Diao's bound and the $2.0996C^{3/4}$ bound}

The following $41$ torus knots satisfy
\[
R_{\mathrm{stick}}(T(q,p))
   >\max\{D(C(T(q,p))),B_{2.0996}(C(T(q,p)))\}.
\]
They are grouped by $p$:
\begin{align*}
p=4:&\quad T(9,4),\\
p=5:&\quad T(9,5),T(11,5),T(12,5),\\
p=6:&\quad T(11,6),T(13,6),\\
p=7:&\quad T(12,7),T(13,7),T(15,7),T(16,7),T(17,7),T(18,7),\\
p=8:&\quad T(13,8),T(15,8),T(17,8),T(19,8),T(21,8),\\
p=9:&\quad T(14,9),T(16,9),T(17,9),T(19,9),T(20,9),T(22,9),T(23,9),\\
p=10:&\quad T(17,10),T(19,10),T(21,10),T(23,10),\\
p=11:&\quad T(17,11),T(18,11),T(19,11),T(20,11),T(21,11),
T(23,11),T(24,11),T(25,11),\\
p=12:&\quad T(19,12),T(23,12),T(25,12),\\
p=13:&\quad T(24,13),T(25,13).
\end{align*}
\subsubsection*{Additional knots beating Diao's bound and the $1.105C^{3/4}$ bound}

The following are the additional knots satisfying
\[
R_{\mathrm{stick}}(T(q,p))
   >\max\{D(C(T(q,p))),B_{1.105}(C(T(q,p)))\},
\]
but not the preceding $2.0996C^{3/4}$ comparison.  To keep the list
compact, for each $p$ we list the corresponding values of $q$:
\[
\mathcal A_{\mathrm{add}}
 =\bigcup_p\{T(q,p):q\in Q_p\}.
\]
Here
\small
\begin{align*}
Q_9={}&\{25\},\\
Q_{10}={}&\{27\},\\
Q_{11}={}&\{26,27,28,29,30,31\},\\
Q_{12}={}&\{29,31\},\\
Q_{13}={}&\{19,20,21,22,23,27,28,29,30,31,32,33,34,35,36,37\},\\
Q_{14}={}&\{23,25,27,29,31,33,37,39,41\},\\
Q_{15}={}&\{22,23,26,28,29,31,32,34,37,38,41,43,44\},\\
Q_{16}={}&\{23,25,27,29,31,33,35,37,39,41,43,45,47\},\\
Q_{17}={}&\{25,26,27,28,29,30,31,32,33,35,36,37,38,39,40,41,42,43,44,45,46,47,48,49,50\},\\
Q_{18}={}&\{29,31,35,37,41,43,47,49,53\},\\
Q_{19}={}&\{27,28,29,30,31,32,33,34,35,36,37,39,40,41,42,43,44,45,46,47,48,49,50,51,52,53,54,55,56\},\\
Q_{20}={}&\{29,31,33,37,39,41,43,47,49,51,53,57,59\},\\
Q_{21}={}&\{31,32,34,37,38,40,41,43,44,46,47,50,52,53,55,58,59,61,62\},\\
Q_{22}={}&\{31,35,37,39,41,43,45,47,49,51,53,57,59,61,63,65\},\\
Q_{23}={}&\{32-45,47,48,49,50,51,52,53,54,55,56,57,58,59,60,61,62,63,64,65,66,67,68\},\\
Q_{24}={}&\{35,37,41,43,47,49,53,55,59,61,65,67,71\},\\
Q_{25}={}&\{36-41,42,43,44,46,47,48,49,51,52,53,54,56,57,58,59,61,62,63,64,66,67,68,69,71-74\},\\
Q_{26}={}&\{37,41,43,45,47,49,51,53,55,57,59,61,63,67,69,71,73\},\\
Q_{27}={}&\{37,38,40,41,43,44,46,47,49,50,52,53,55,56,58,59,61,62,64,65,67,68,70,71\},\\
Q_{28}={}&\{39,41,43,45,47,51,53,55,57,59,61,65,67\},\\
Q_{29}={}&\{40,41,42,43,44,45,46,47,48,49,50,51,52,53,54,55,56,57,59,60,61,62,63,64,65,66\},\\
Q_{30}={}&\{41,43,47,49,53,59,61\},\\
Q_{31}={}&\{42,43,44,45,46,47,48,49,50,51,52,53,54,55,56,57,58,59,60,61\},\\
Q_{32}={}&\{45,47,49,51,53,55,57,59\},\\
Q_{33}={}&\{46,47,49,50,52,53,56\},\\
Q_{34}={}&\{47,49,53,55\},\\
Q_{35}={}&\{48,51,52,53,54\},\\
Q_{36}={}&\{49\},\\
Q_{37}={}&\{50,51\}.
\end{align*}
\normalsize
There are $387$ knots in this additional list.  Thus, within the stated
crossing-number range and stick-number theorem, the $2.0996$ comparison
contains $41$ knots, while the weaker $1.105$ comparison contains
$41+387=428$ knots in total.

\subsection{A conjectural stick-based lower bound for \(T(P,2)\)}

For the alternating torus knots \(T(P,2)\), with \(P\) odd, the
conjectured stick-number formula is
\[
s\bigl(T(P,2)\bigr)
=
\left\lfloor \frac{2(P+1)}{3}\right\rfloor+4.
\tag{2}
\]
The available polygonal constructions attain the right-hand side, but
the corresponding general lower bound has not been proved.

\textbf{Author's note: Cantarella et al. \cite{CantarellaRechnitzerSchumacherShonkwiler2026} state that this could not be improved upon by their algorithm for $P$ from 3 to 31.}

Substituting the conjectured value (2) into (1) gives the conjectural
ropelength lower bound
\[
\boxed{
\operatorname{Rop}\bigl(T(P,2)\bigr)
\geq
\pi\left(
\left\lfloor \frac{2(P+1)}{3}\right\rfloor+3
\right).
}
\tag{3}
\]

Since \(P\) is the crossing number of \(T(P,2)\),
\[
C\bigl(T(P,2)\bigr)=P,
\]
and therefore
\[
\operatorname{Rop}\bigl(T(P,2)\bigr)
\geq
\frac{2\pi}{3}P+O(1)
=
\frac{2\pi}{3}C+O(1).
\tag{4}
\]
Thus the conjectured asymptotic slope is
\[
\boxed{
\liminf_{P\to\infty}
\frac{\operatorname{Rop}(T(P,2))}{P}
\geq
\frac{2\pi}{3}
\approx 2.094395.
}
\tag{5}
\]

Removing the floor gives the slightly weaker explicit inequality
\[
\operatorname{Rop}\bigl(T(P,2)\bigr)
\geq
\frac{2\pi}{3}P+\frac{8\pi}{3}.
\tag{6}
\]

The ropelength implication is rigorous conditional on the stick-number
inequality in (2).  Hence the status of (3)--(6) is conjectural only
because the general minimality of the proposed stick number remains
unproved.

\textbf{Author's note: the LLM claimed for a good part of our dialogue this construction was an exact lower bound for alternating torus knots, which would have made this a much stronger result. This is one of many warnings that LLM outputs should still be scrutinized.}

\subsection{Pretzel knots and stick-based ropelength bounds}

Consider the alternating pretzel family
\[
P_m=P(2,\underbrace{3,3,\ldots,3}_{m\text{ entries}}),
\qquad m\geq 2.
\]
Its standard diagram is reduced and alternating, so its crossing number is
\[
c(P_m)=3m+2.
\]
The bridge number is
\[
b(P_m)=m+1.
\]

For a unit-thickness knot, the bridge-number estimate gives
\[
\operatorname{Rop}(K)\geq 2\pi b(K),
\]
and hence
\[
\operatorname{Rop}(P_m)\geq 2\pi(m+1).
\]

The bridge and superbridge indices satisfy
\[
b(K)<\operatorname{sb}(K)
\leq \frac{1}{2}\operatorname{stick}(K).
\]
Since these indices are integers,
\[
\operatorname{stick}(P_m)\geq 2b(P_m)+2=2m+4.
\]
The standard stick--ropelength inequality,
\[
\operatorname{stick}(K)
\leq
\left\lfloor\frac{\operatorname{Rop}(K)}{\pi}\right\rfloor+1,
\]
therefore implies
\[
\operatorname{Rop}(P_m)
\geq
\pi\bigl(\operatorname{stick}(P_m)-1\bigr)
\geq \pi(2m+3).
\]

In terms of crossing number,
\[
\boxed{
\operatorname{Rop}(P_m)
\geq
\frac{2\pi}{3}\,c(P_m)+\frac{5\pi}{3}.
}
\]
This improves the direct bridge-number estimate
\[
\operatorname{Rop}(P_m)
\geq
\frac{2\pi}{3}\,c(P_m)+\frac{4\pi}{3}
\]
by exactly \(\pi\). The stick bound is therefore the stronger of the two bridge-based estimates for this family.



\subsection{A non-pretzel Montesinos family}

A Montesinos knot is obtained by taking the numerator closure of a cyclic sum
of rational tangles.  Consider the family
\[
K_n=
M\left(
\frac{1}{2},
\underbrace{\frac{2}{5},\frac{2}{5},\ldots,\frac{2}{5}}_{n-1\text{ tangles}}
\right),
\qquad n\geq 3.
\]
Because exactly one denominator is even, this is a knot rather than a link.
Since the tangle \(1/2\) has a reduced alternating diagram with two
crossings and \(2/5\) has one with four crossings, the standard
all-positive Montesinos diagram has
\[
c(K_n)=2+4(n-1)=4n-2
\]
crossings.  The diagram is reduced and alternating, so this is the
crossing number.

For a normalized Montesinos knot with \(n\) rational tangles, the bridge
number is \(n\).  Hence
\[
b(K_n)=n.
\]
The bridge, superbridge, and stick indices satisfy
\[
b(K)<\operatorname{sb}(K)
\leq \frac{1}{2}\operatorname{stick}(K).
\]
Since these indices are integral, it follows that
\[
\operatorname{stick}(K_n)\geq 2b(K_n)+2=2n+2.
\]

For a unit-radius rope, the stick--ropelength inequality is
\[
\operatorname{stick}(K)
\leq
\left\lfloor\frac{\operatorname{Rop}(K)}{\pi}\right\rfloor+1.
\]
Consequently,
\[
\operatorname{Rop}(K_n)
\geq
\pi\bigl(\operatorname{stick}(K_n)-1\bigr)
\geq \pi(2n+1).
\]
Using \(c(K_n)=4n-2\), this becomes
\[
\boxed{
\operatorname{Rop}(K_n)
\geq
\frac{\pi}{2}\,c(K_n)+2\pi.
}
\]

For comparison, the direct bridge-number estimate gives only
\[
\operatorname{Rop}(K_n)\geq 2\pi b(K_n)=2\pi n
=\frac{\pi}{2}\,c(K_n)+\pi.
\]
Thus the stick argument improves the bridge estimate by exactly \(\pi\).

The stick-based bound first exceeds Diao's crossing-number estimate at
\(n=7\), for which \(c(K_7)=26\):
\[
\operatorname{Rop}(K_7)\geq 15\pi\approx47.124,
\]
whereas Diao's bound is approximately \(45.843\).  Since the
Montesinos estimate grows linearly in \(c(K_n)\), it remains stronger
for all \(n\geq 7\).

This argument is rigorous provided the standard results on the crossing
and bridge numbers of normalized Montesinos knots are invoked.  It does
not require the exact stick number of \(K_n\); only the lower bound
\(\operatorname{stick}(K_n)\geq2n+2\) is used.

\textbf{Author's note: in both these cases, the linear slope based on the bridge number could have already been deduced from published literature \cite{DiaoErnstThistlethwaite2003, CantarellaRechnitzerSchumacherShonkwiler2026}, but had not been explicitly stated. The stick number bound simply adds $\pi$ to the lower bound, but this is still the first time the strong crossing-dependence of these knots has been stated.}

\subsection{Asymmetric-helix construction of $P(2,3^N)$}
\label{subsec:scripted-asymmetric-pretzels}

We  can construct a polygonal approximation to the pretzel knot
with twist parameters $P(2,3,\ldots,3)$, containing $N$ entries equal to $3$.
Write $M=N+1$ for the number of twist modules, and index them by
$i=0,\ldots,N$.  The first module has $p_0=2$ half-twists and the remaining
modules have $p_i=3$.  The modules are placed at equally spaced points on a
circle in the $xy$-plane.  With
\[
 D=0.92\left(2(8.05)\sin\frac{\pi}{7}\right),\qquad
 \rho_N=\frac{D}{2\sin(\pi/M)},\qquad
 \alpha_i=\frac{2\pi i}{M},
\]
the center and local frame of module $i$ are
\[
 \mathbf c_i=\rho_N(\cos\alpha_i,\sin\alpha_i,0),\quad
 \mathbf u_i=(-\sin\alpha_i,\cos\alpha_i,0),\quad
 \mathbf v_i=(\cos\alpha_i,\sin\alpha_i,0).
\]
Here $\mathbf u_i$ is tangent to the circle of module centers and $\mathbf v_i$
is its radial direction.

\textbf{Author's note: I am avoiding revealing the exact parameters of the helix because I am aware of some human authors trying to publish a paper about them. These parameters were deduced independently by ChatGPT, and the human authors have already found a better solution than that. I imagine norms regarding this choice are evolving rapidly. Here, I will use the parameters from Kim et al. \cite{kim2024efficiency} instead of the slightly better ones that I hope my colleagues are able to publish soon.}

Each module is represented by two asymmetric helical arcs.  Let $s=+1$ denote
the small-radius strand and $s=-1$ the large-radius strand, and set

\[
 a_{+1}=r_{\rm in}=0.142,\qquad
 a_{-1}=-R_{\rm out}=1.9498,\qquad H=0.3226.
\]
For $0\leq t\leq1$, the centerline of a strand is
\[
 \boldsymbol\gamma_{i,s}(t)=\mathbf c_i+
 a_s\bigl(\cos\theta_{i,s}(t)\,\mathbf u_i+
 \sin\theta_{i,s}(t)\,\mathbf v_i\bigr)+z_i(t)\,\mathbf e_z,
 \qquad \theta_{i,s}(t)=\phi_{i,s}+\omega_i t.
\]
The phase and axial coordinate are assigned as follows:
\[
 (\omega_i,\phi_{i,s})=
 \begin{cases}
 (7\pi/4,\,0), & i=0,\ s=+1,\\
 (7\pi/4,\,0.72), & i=0,\ s=-1,\\
 (5\pi/2,\,\pi/4), & 1\leq i\leq N,\ s\in\{-1,+1\},
 \end{cases}
 \qquad
 z_i(t)=
 \begin{cases}
 -\pi H+(7\pi H/4)t, & i=0,\\
 3\pi H(t-1/2), & 1\leq i\leq N.
 \end{cases}
\]
Thus the $3$-modules use the shortened phase interval
$\pi/4\leq\theta\leq11\pi/4$, while the $2$-module uses a $7\pi/4$ phase
sweep and a separate offset for its large-radius strand.  The two radii and
opposite signs of $a_s$ give the asymmetric pair of helices used in the scripts.

The helix endpoints are joined according to the pretzel endpoint pairing.  A
strand entering module $i$ at the bottom side $s$ leaves at the top side
$s(-1)^{p_i}$.  For adjacent modules, including the cyclic pair $N$ and $0$,
the scripts connect the top endpoint of side $+1$ of module $i$ to the top
endpoint of side $-1$ of module $i+1$, and make the analogous connection at the
bottom.  This pairing produces one closed component for the instances accepted
by the scripts.

Each connector is a cubic Bezier arc with endpoint tangents taken from the
helices.  Its two handle lengths are searched to reduce connector length subject
to a sampled minimum radius of curvature of $1.05$.  If that search finds no
feasible cubic, the scripts try a quartic bulge of the form
$16t^2(1-t)^2\mathbf q$ in several directions.  The connectors are optimized
individually; the scripts do not impose a global minimum-distance constraint
between all parts of the knot.

The helix and connector arcs are first sampled densely and concatenated into
one closed curve.  The cumulative polygonal arclength is then used to resample
the curve at $n$ equally spaced locations.  The output is an $n\times3$ array
of Cartesian coordinates, written as tab-separated values.  The associated
pretzel diagram has $2+3N$ crossings; for example, $N=10$ gives the standard
$P(2,3^{10})$ diagram with $32$ crossings.  The script's reported polygonal
length is the length of the sampled output; the coordinates are not
automatically rescaled to a prescribed global thickness.

For $N=10$, the method place $11$ modules on a circle of radius
$\rho_{10}\approx11.405650$, with adjacent module centers separated by
$D\approx6.426686$.










\section{Optimal parameters for the standard torus knot parameterization}

\subsection{Radial and vertical scaling in the standard torus family}
\label{sec:standard-torus-scalings}

We normalize the tube radius (and hence the minor radius of the
parameterizing torus) to one.  For coprime integers \(P>Q\geq 2\), consider
the standard torus parameterization
\begin{equation}
 \gamma_{P,Q}^{(\rho,\lambda)}(t)
 =
 \bigl((\rho+\cos Pt)\cos Qt,\,
       (\rho+\cos Pt)\sin Qt,\,
       \lambda\sin Pt\bigr),
 \qquad 0\leq t\leq 2\pi .
\label{eq:standard-torus}
\end{equation}
Here \(\rho=R/r\) is the major-to-minor radius ratio and \(\lambda\) scales
the vertical coordinate.  The choice \(\lambda=1\) gives a circular
cross-section.  The curve is a \(T(P,Q)\) knot when \(\gcd(P,Q)=1\); when
\(\gcd(P,Q)>1\), it is the corresponding torus link.

\subsection{Length and the radial ansatz}

Put \(u=Pt\).  In cylindrical coordinates, the radial coordinate is
\(a(u)=\rho+\cos u\), the azimuthal angle is \(Qt\), and the vertical
coordinate is \(\lambda\sin u\).  Consequently,
\begin{equation}
 \left|\gamma'(t)\right|^2
 =
 P^2\bigl(\sin^2 u+\lambda^2\cos^2 u\bigr)
 +Q^2(\rho+\cos u)^2 .
\label{eq:speed}
\end{equation}
Since the integrand is \(2\pi\)-periodic in \(u\),
\begin{equation}
 L(P,Q;\rho,\lambda)
 =
 \int_0^{2\pi}
 \sqrt{
 P^2\bigl(\sin^2 u+\lambda^2\cos^2 u\bigr)
 +Q^2(\rho+\cos u)^2
 }\,du .
\label{eq:length-integral}
\end{equation}

The simple radial rule used below is
\begin{equation}
 \boxed{\rho_{\rm rad}=1+\frac{P}{Q}.}
\label{eq:radial-rule}
\end{equation}
Geometrically, the inner radius of the torus is then
\(\rho_{\rm rad}-1=P/Q\).  Thus the radius available to the rapid
longitudinal motion is balanced against the ratio of the meridional and
longitudinal frequencies.  Equation~\eqref{eq:radial-rule} is a useful
explicit geometric ansatz; it is not a theorem asserting global optimality
among all embeddings of the knot type.

There is also a direct local interpretation of this rule.  At the inner
side of the torus, \(u=\pi\), the meridional contribution to the speed is
\(P\), whereas the longitudinal contribution is \(Q(\rho-1)\).  Equating
these two characteristic speeds gives
\[
 P=Q(\rho-1),
\]
which is precisely \eqref{eq:radial-rule}.

\subsection{Geometric vertical-scale rules}

Two convenient one-parameter vertical corrections are
\begin{equation}
 \lambda_{\sqrt{\vphantom{a}}}
 =
 \sqrt{1+\left(\frac{Q}{P}\right)^2}
 =
 \frac{\sqrt{P^2+Q^2}}{P},
\qquad
 \lambda_{\csc}
 =
 \csc\left(\frac{\pi Q}{P}\right).
\label{eq:vertical-rules}
\end{equation}
The square-root rule can be obtained by replacing the meridional frequency
\(P\) in the vertical direction by the Euclidean combined frequency
\(\sqrt{P^2+Q^2}\).  In other words, it is the normalization determined by
\[
 P^2\lambda_{\sqrt{\vphantom{a}}}^{\,2}=P^2+Q^2.
\]
This makes the vertical speed at \(u=0\) reflect both angular frequencies,
while retaining the radial balance \(\rho-1=P/Q\).

The cosecant rule is a chord-compensation prescription.  A half-turn in
the meridional variable, \(\Delta t=\pi/P\), advances the longitudinal
angle by \(\pi Q/P\).  If the transverse scale is reduced by the sine of
this angular advance, restoring unit transverse separation requires
\[
 \lambda_{\csc}\sin\left(\frac{\pi Q}{P}\right)=1,
\]
which gives the second formula in \eqref{eq:vertical-rules}.  Both rules
are geometric ansatzes rather than consequences of a global variational
theorem; their quality is determined by the resulting thickness in
\eqref{eq:thickness}.

\textbf{Author's note: plots of the ropelengths of these configurations (when selecting the optimal $\lambda$) against crossing number can be seen in Figure 1, and a plot of the optimal $\lambda$ and its effect against $P/Q$ can be seen in Figure 4. The apparent groupings of the points into distinct lines correspond to different fixed $Q$ with varying $P$. These slopes are discussed subsequently, and are known to always exceed $4Q$.}

\subsection{Appell \(F_1\) formula for the length}

For fixed \((P,Q,\rho,\lambda)\), define
\begin{equation}
\begin{split}
 D&=P^2\lambda^2+Q^2(\rho+1)^2,\\
 A&=P^2(\lambda^2-1)+Q^2(\rho+1),\\
 B&=P^2(\lambda^2-1)+Q^2 .
\end{split}
\label{eq:DAB}
\end{equation}
On setting \(u=2v\), the expression under the square root in
\eqref{eq:length-integral} factors as
\begin{equation}
 P^2\bigl(\sin^2 u+\lambda^2\cos^2 u\bigr)
 +Q^2(\rho+\cos u)^2
 =
 D(1-x_+\sin^2v)(1-x_-\sin^2v),
\label{eq:factorization}
\end{equation}
where
\begin{equation}
 x_\pm
 =
 \frac{2}{D}\left(A\pm\sqrt{A^2-BD}\right).
\label{eq:xpm}
\end{equation}
The square root in \eqref{eq:xpm} may be imaginary; the two conjugate
terms still give a real positive integrand.

The Euler integral for Appell's hypergeometric function \(F_1\) gives
\begin{equation}
 \int_0^{\pi/2}
 (1-x\sin^2v)^{1/2}(1-y\sin^2v)^{1/2}\,dv
 =
 \frac{\pi}{2}
 F_1\!\left(\frac12;-\frac12,-\frac12;1;x,y\right).
\label{eq:appell-integral}
\end{equation}
Therefore the length of the standard torus curve is
\begin{equation}
 \boxed{
 L(P,Q;\rho,\lambda)
 =
 2\pi\sqrt D\,
 F_1\!\left(\frac12;-\frac12,-\frac12;1;x_+,x_-\right).
 }
\label{eq:appell-length}
\end{equation}
For radial scaling alone, put \(\lambda=1\) and
\(\rho=1+P/Q\) in \eqref{eq:appell-length}.  In this case
\begin{equation}
 D=P^2+Q^2(\rho+1)^2,\qquad
 A=Q^2(\rho+1),\qquad B=Q^2 .
\label{eq:radial-DAB}
\end{equation}

\subsection{Thickness and Appell-based ropelength formulae}

For a \(C^2\) embedded curve, write \(\tau\) for its thickness.  Equivalently,
\begin{equation}
 \tau(\gamma)
 =
 \min\left\{
 \inf_t\frac{|\gamma'(t)|^3}
 {|\,\gamma'(t)\mathbin{\times}\gamma''(t)\,|},
 \;
 \frac12\inf_{\substack{s\ne t\\
 (\gamma(t)-\gamma(s))\cdot\gamma'(t)=0\\
 (\gamma(t)-\gamma(s))\cdot\gamma'(s)=0}}
 |\gamma(t)-\gamma(s)|
 \right\}.
\label{eq:thickness}
\end{equation}
The first term is the minimum radius of curvature and the second is the
doubly-critical self-distance.  Thus the exact ropelength of this
parameterized curve is
\begin{equation}
 \boxed{
 \mathcal R(P,Q;\rho,\lambda)
 =
 \frac{
 2\pi\sqrt D\,
 F_1\!\left(\frac12;-\frac12,-\frac12;1;x_+,x_-\right)
 }{\tau(P,Q;\rho,\lambda)} .
 }
\label{eq:appell-ropelength}
\end{equation}
For reference, if \(a=\rho+\cos u\), the speed and curvature numerator
can be written
\begin{equation}
 v^2=P^2(\sin^2u+\lambda^2\cos^2u)+Q^2a^2,
\end{equation}
\begin{equation}
\begin{split}
 N^2=|\gamma'\mathbin{\times}\gamma''|^2
={}&\lambda^2P^4Q^2\sin^2u(\cos u-\rho)^2\\
&+\lambda^2P^2\bigl(P^2+Q^2a\cos u\bigr)^2\\
&+Q^2\bigl(2P^2\sin^2u+P^2a\cos u+Q^2a^2\bigr)^2 ,
\end{split}
\label{eq:curvature-numerator}
\end{equation}
so that
\begin{equation}
 \tau_{\rm curv}=\min_{0\leq u\leq 2\pi}\frac{v^3}{N}.
\label{eq:curvature-thickness}
\end{equation}

The three explicit candidates considered here are
\begin{align}
 \mathcal R_{\rm rad}(P,Q)
 &=\mathcal R\left(P,Q;1+\frac{P}{Q},1\right), \label{eq:Rrad}\\
 \mathcal R_{\sqrt{\vphantom{a}}}(P,Q)
 &=\mathcal R\left(P,Q;1+\frac{P}{Q},
          \sqrt{1+\left(\frac{Q}{P}\right)^2}\right), \label{eq:Rsqrt}\\
 \mathcal R_{\csc}(P,Q)
 &=\mathcal R\left(P,Q;1+\frac{P}{Q},
          \csc\left(\frac{\pi Q}{P}\right)\right). \label{eq:Rcsc}
\end{align}
The last two are vertical-scaling ansatzes.  They are explicit constructions
whose admissibility must still be checked through \eqref{eq:thickness};
they are not claims of global ropelength minimization.

\subsection{Fixed-\(Q\) slopes as \(P\) grows}

Consider the radial family with \(Q\) fixed and \(P\to\infty\).  The
length has the leading behavior
\begin{equation}
 \frac{L(P,Q;1+P/Q,1)}{P}
 \longrightarrow 2\pi\sqrt2 .
\label{eq:fixed-Q-length}
\end{equation}
The thickness is controlled by a near-repeat doubly-critical pair.  A
convenient limiting description uses
\(\theta_Q=\pi/Q\) and the root \(\beta_Q\) near
\(-\theta_Q/2\) of
\begin{equation}
 \beta_Q+
 \sin(\theta_Q+\beta_Q)\cos(\theta_Q+\beta_Q)=0.
\label{eq:beta-equation}
\end{equation}
The corresponding limiting full-distance coefficient is
\begin{equation}
 d_Q=
 \sqrt{\beta_Q^2+\sin^2(\theta_Q+\beta_Q)},
\label{eq:dQ}
\end{equation}
and the limiting thickness is \(d_Q/(2P)\).  Hence
\begin{equation}
 \boxed{
 \mathcal R_{\rm rad}(P,Q)
 \sim m_QP,\qquad
 m_Q=\frac{2\pi\sqrt2}{d_Q}.
 }
\label{eq:slope}
\end{equation}
Numerically, this gives
\begin{center}
\begin{tabular}{c|cccccccc}
\(Q\) & 3&4&5&6&7&8&9&10\\ \hline
\(m_Q\) &12.3396&16.2297&20.1762&24.1436&
28.1216&32.1055&36.0933&40.0836
\end{tabular}
\end{center}
The values are limiting numerical evaluations of
\eqref{eq:beta-equation}--\eqref{eq:slope}, not universal lower bounds for
arbitrary \(T(P,Q)\) embeddings.

\subsection{Large-\(Q\) series for the fixed-\(Q\) slope}

Let \(\theta=\pi/Q\).  Expanding the root in
\eqref{eq:beta-equation} gives
\begin{equation}
 \beta_Q
 =
 -\frac{\theta}{2}
 +\frac{\theta^3}{24}
 +\frac{\theta^5}{120}
 +O(\theta^7).
\label{eq:beta-series}
\end{equation}
Substitution into \eqref{eq:dQ} yields
\begin{equation}
 d_Q
 =
 \frac{\theta}{\sqrt2}
 \left(
 1-\frac{\theta^2}{48}
 -\frac{23\theta^4}{7680}
 +O(\theta^6)
 \right).
\label{eq:d-series}
\end{equation}
Inverting this expansion in \eqref{eq:slope} gives
\begin{equation}
 \boxed{
 m_Q
 =
 4Q+\frac{\pi^2}{12Q}
 +\frac{79\pi^4}{5760Q^3}
 +O(Q^{-5}).
 }
\label{eq:slope-series}
\end{equation}
Thus the commonly useful approximation
\begin{equation}
 m_Q\approx 4Q+\frac{\pi^2}{12Q}
\label{eq:two-term-slope}
\end{equation}
is the first two terms of the large-\(Q\) expansion.  The next displayed
term is positive, so the two-term expression is asymptotically slightly
below the limiting slope; it should not be described as an overestimate
without an additional finite-\(Q\) argument.

\subsection{The limiting \(T(P,P-1)\) vertical correction}

Now let \(Q=P-1\) and let \(P\to\infty\).  The radial ratio tends to
\(\rho=2\).  For a fixed vertical scale \(\lambda\), define
\begin{equation}
 I(\lambda)
 =
 \int_0^{2\pi}
 \sqrt{
 \sin^2u+\lambda^2\cos^2u+(2+\cos u)^2
 }\,du .
\label{eq:I-lambda}
\end{equation}
The length satisfies
\begin{equation}
 L(P,P-1;2,\lambda)
 \sim P I(\lambda).
\label{eq:P-Pminus1-length}
\end{equation}

The relevant near-repeat pair has parameter separation
\(\Delta t=2\pi/(P-1)+O(P^{-2})\).  Expanding the squared distance to
the first nontrivial order in \(P^{-1}\), and writing the phase of the
pair as \(u\), gives
\begin{equation}
 D_\lambda^2(u,a)
 =
 a^2\bigl(\sin^2u+\lambda^2\cos^2u\bigr)
 +(2\pi+a)^2(2+\cos u)^2 ,
\label{eq:PP-distance}
\end{equation}
where \(a\) is the bounded first-order correction to the near-repeat
separation.  For \(\lambda=1\) and \(\lambda=\sqrt2\), the minimum in the
phase occurs at the inner point \(u=\pi\).  Minimizing there with respect
to \(a\) gives the limiting full-distance coefficients
\begin{equation}
 d_1=\frac{2\pi}{\sqrt2}=\pi\sqrt2,
\qquad
 d_{\sqrt2}=\frac{2\pi\sqrt2}{\sqrt3}.
\label{eq:d-P-Pminus1}
\end{equation}
Since thickness is half the full doubly-critical distance,
\(\tau_\lambda\sim d_\lambda/(2P)\).  Therefore the ratio of the two
ropelengths is
\begin{equation}
 \frac{\mathcal R_{\sqrt{\vphantom{a}}}(P,P-1)}
      {\mathcal R_{\rm rad}(P,P-1)}
 \longrightarrow
 \frac{I(\sqrt2)}{I(1)}
 \frac{d_1}{d_{\sqrt2}}
 =
 \frac{I(\sqrt2)}{I(1)}\frac{\sqrt3}{2}.
\label{eq:PPminus1-ratio}
\end{equation}
Numerically,
\begin{equation}
 I(1)=14.2073964944,\qquad
 I(\sqrt2)=14.9422807277,
\end{equation}
so that
\begin{equation}
 \frac{I(\sqrt2)}{I(1)}=1.0517254680,
\qquad
 \frac{\mathcal R_{\sqrt{\vphantom{a}}}}
      {\mathcal R_{\rm rad}}
 \longrightarrow 0.9108209731.
\label{eq:PPminus1-numerics}
\end{equation}
The limiting reduction is consequently
\begin{equation}
 \boxed{
 1-0.9108209731
 =0.0891790269
 \approx 8.92\%\approx 9\%.
 }
\label{eq:nine-percent}
\end{equation}
This is a comparison of two explicit standard-torus constructions.  It is
not a statement that the vertically scaled curve is a global ropelength
minimizer.  Finally, the cosecant ansatz behaves differently in this
limit:
\(\csc(\pi Q/P)\sim P/\pi\) when \(Q=P-1\), so it does not have the same
finite limiting correction as the square-root scaling.

\subsection{A general upper bound for the ropelength of torus knots}
\label{sec:two-denominator-upper-bounds}

Assume \(P>Q\geq 2\) and \(\gcd(P,Q)=1\).  The same formulas apply
to the corresponding torus link when \(\gcd(P,Q)>1\).  Normalize the
minor radius of the parameterizing torus to one and consider
\begin{equation}
 \gamma_{P,Q}(t)
 =
 \bigl((\rho+\cos Pt)\cos Qt,\,
       (\rho+\cos Pt)\sin Qt,\,
       \sin Pt\bigr),
 \qquad 0\leq t\leq 2\pi ,
\end{equation}
with the radial choice
\begin{equation}
 \rho=1+\frac{P}{Q}.
\end{equation}

The speed is
\begin{equation}
 |\gamma_{P,Q}'(t)|^2
 =
 P^2+Q^2(\rho+\cos Pt)^2.
\end{equation}
Since the square-root function is concave,
\begin{align}
 L(P,Q)
 &\leq
 2\pi\sqrt{P^2+Q^2\left(\rho^2+\frac12\right)} \\
 &=
 2\pi\sqrt{2P^2+2PQ+\frac32Q^2}.
 \label{eq:torus-length-upper-bound}
\end{align}

Define
\begin{equation}
 B_{P,Q}
 =
 \sqrt{
 2P^4+2P^3Q+9P^2Q^2+4PQ^3+4Q^4
 },
\end{equation}
and
\begin{equation}
 V_{P,Q}
 =
 \sqrt{2P^2+4PQ+4Q^2}.
\end{equation}

The minimum speed satisfies
\[
 |\gamma_{P,Q}'(t)|^2\geq 2P^2,
\]
while a direct estimate of the acceleration gives
\[
 |\gamma_{P,Q}''(t)|\leq B_{P,Q}.
\]
Consequently, the curvature obeys
\[
 \kappa(t)\leq \frac{B_{P,Q}}{2P^2},
\]
so the curvature radius is bounded below by
\begin{equation}
 D_{\kappa}(P,Q)
 =
 \frac{2P^2}{B_{P,Q}}.
\end{equation}

For two parameter values separated by \(h\), the chord distance satisfies
\begin{equation}
 \left|\gamma_{P,Q}(t+h)-\gamma_{P,Q}(t)\right|^2
 \geq
 4\left[
 \sin^2\left(\frac{Ph}{2}\right)
 +
 \left(\frac{P}{Q}\right)^2
 \sin^2\left(\frac{Qh}{2}\right)
 \right].
 \label{eq:torus-chord-lower-bound}
\end{equation}

Using the coprimality of \(P\) and \(Q\), together with the minimum
parameter separation of a doubly-critical pair, gives the conservative
nonlocal thickness bound
\begin{equation}
 D_{\mathrm{nl}}(P,Q)
 =
 \frac{\sqrt{2}}{Q}
 \quad\text{or}\quad
 D_{\mathrm{sep}}(P,Q)
 =
 \frac{\sqrt{2}\,P^3}{B_{P,Q}V_{P,Q}},
 \label{eq:two-nonlocal-denominators}
\end{equation}
with the smaller of these two quantities controlling the thickness.

The curvature term is never the active term in this estimate, because
\begin{equation}
 \frac{D_{\mathrm{sep}}(P,Q)}{D_{\kappa}(P,Q)}
 =
 \frac{P}{\sqrt{2}\,V_{P,Q}}
 <\frac12.
\end{equation}

Therefore, only two denominator possibilities remain.

Combining \eqref{eq:torus-length-upper-bound} and
\eqref{eq:two-nonlocal-denominators}, one obtains
\begin{equation}
\boxed{
 U(P,Q)
 =
 \frac{
 2\pi\sqrt{2P^2+2PQ+\frac32Q^2}
 }{
 \min\left\{
 \dfrac{\sqrt{2}}{Q},
 \;
 \dfrac{\sqrt{2}\,P^3}
 {
 \sqrt{2P^4+2P^3Q+9P^2Q^2+4PQ^3+4Q^4}
 \sqrt{2P^2+4PQ+4Q^2}
 }
 \right\}
 }.
}
\label{eq:UPQ-two-denominator}
\end{equation}

The first denominator is active when
\begin{equation}
 B_{P,Q}V_{P,Q}\leq QP^3,
\end{equation}
and the second is active when
\begin{equation}
 B_{P,Q}V_{P,Q}\geq QP^3.
\end{equation}

\textbf{Author's note: I thought it would be useful to have an upper bound as a function of crossing number at a given ratio of P/Q, given my comparison of $T(Q,Q-1)$, $T(Q,Q)$, $T(2Q,Q)$ etc. However, this should be treated as approximate for knots as the co-primality of $P$ and $Q$ makes it unlikely for the ratio to stay constant as $C$ increases.}
Let
\begin{equation}
 C=P(Q-1),
 \qquad
 r=\frac{P}{Q}.
\end{equation}
Then
\begin{equation}
 Q=\frac{1+\sqrt{1+4C/r}}{2}.
\end{equation}

Define
\begin{equation}
 b(r)=\sqrt{2r^4+2r^3+9r^2+4r+4},
 \qquad
 v(r)=\sqrt{2r^2+4r+4}.
\end{equation}

Substitution into \eqref{eq:UPQ-two-denominator} gives
\begin{equation}
\boxed{
 U(C,r)
 =
 \frac{
 \pi\left(1+\sqrt{1+\frac{4C}{r}}\right)
 \sqrt{2r^2+2r+\frac32}
 }{
 \min\left\{
 \dfrac{2\sqrt{2}}
 {1+\sqrt{1+\frac{4C}{r}}},
 \;
 \dfrac{\sqrt{2}\,r^3}{b(r)v(r)}
 \right\}
 }.
}
\label{eq:UCr-two-denominator}
\end{equation}

For fixed \(r\), the first denominator tends to zero like \(C^{-1/2}\),
whereas the second denominator is constant.  Hence the first denominator
is eventually active and
\begin{align}
 U(C,r)
 &=
 \frac{
 \pi\left(1+\sqrt{1+\frac{4C}{r}}\right)^2
 \sqrt{2r^2+2r+\frac32}
 }{2\sqrt{2}} \\
 &=
 \pi\sqrt{4+\frac4r+\frac{3}{r^2}}\,C
 +
 \sqrt{2}\pi
 \sqrt{2r^2+2r+\frac32}
 \sqrt{\frac{C}{r}}
 +O(1).
 \label{eq:fixed-r-limit}
\end{align}
In particular,
\begin{equation}
\boxed{
 U(C,r)
 \sim
 \pi\sqrt{4+\frac4r+\frac{3}{r^2}}\,C
 \qquad (C\to\infty,\ r\ \mathrm{fixed}).
}
\end{equation}

For fixed \(Q\), one has
\begin{equation}
 B_{P,Q}\sim\sqrt{2}\,P^2,
 \qquad
 V_{P,Q}\sim\sqrt{2}\,P.
\end{equation}
Thus
\[
 D_{\mathrm{sep}}(P,Q)\longrightarrow\frac1{\sqrt2},
 \qquad
 D_{\mathrm{ang}}(P,Q)=\frac{\sqrt2}{Q}.
\]

For \(Q\geq3\), the angular denominator is eventually smaller, and
\begin{equation}
\boxed{
 U(P,Q)
 =
 2\pi QP+\pi Q^2+O(P^{-1}).
}
\label{eq:fixed-Q-upper-bound}
\end{equation}

For \(Q=2\), the two denominators have the same limiting value, but the
finite-separation denominator is slightly smaller for finite \(P\).  In
this case,
\begin{equation}
\boxed{
 U(P,2)
 =
 4\pi P+16\pi+O(P^{-1}).
}
\end{equation}

Since \(C=P(Q-1)\), the corresponding crossing-number forms are
\begin{equation}
\boxed{
 U(C,Q)
 =
 \frac{2\pi Q}{Q-1}\,C+\pi Q^2+O(C^{-1}),
 \qquad Q\geq3,
}
\end{equation}
and
\begin{equation}
\boxed{
 U(C,2)
 =
 4\pi C+16\pi+O(C^{-1}).
}
\end{equation}

These estimates are rigorous upper bounds for the explicitly parameterized
torus family under the stated thickness estimates.  They are not claims
that the standard torus parameterization is globally ropelength-minimizing.

\section{Improvements to $T(Q,Q)$ torus link constructions and bounds}
\textbf{Author's note: I recommend reading my paper for context \cite{klotz2026tight}. It concerns creating tight $T(Q,Q)$ links out of two sets of concentric toroidal helices, Hopf linked with each other.}

\subsection{Upper bound determination for optimal $T(pQ,Q)$ for $p>1$}
\label{sec:generalized-double-donut}

Using Klotz's construction, the $p=1$ value $11.64$ is obtained by
maximizing the number of helices in each shell, integrating the shell
lengths, doubling the construction, and dividing by the crossing number of
the double donut.  The double-donut crossing count is
\begin{equation}
 C_{\mathrm{dbl}}
 =2pQ(Q-1)+2Q^2
 \sim 2(p+1)Q^2.
 \label{eq:double-donut-crossing-supplied}
\end{equation}

\subsubsection{Shell capacity}

Let $T$ be the number of shells, with shell radii
\begin{equation}
 r_i=2i,\qquad 1\leq i\leq T,
 \label{eq:general-shell-radius}
\end{equation}
so that the outer radius is $r_o=2T$.  Let the donut-hole radius be
\begin{equation}
 h=2\lambda T,\qquad \lambda\geq 1.
 \label{eq:general-hole-radius}
\end{equation}
The case $\lambda=1$ is the smallest hole that allows two identical donuts to
pass through one another asymptotically.

For a $p$-fold helix, the small-angle no-overlap equation becomes
\begin{equation}
 d^2(\theta)
 \approx r_i^2\left(p\theta+\frac{2\pi}{N_i}\right)^2+h^2\theta^2.
 \label{eq:p-fold-rectangle-distance}
\end{equation}
Requiring the minimum distance to equal $2$ gives
\begin{equation}
 N_i
 \sim \frac{\pi h r_i}{\sqrt{h^2+p^2r_i^2}}
 =\frac{2\pi\lambda iT}
 {\sqrt{\lambda^2T^2+p^2i^2}}.
 \label{eq:p-fold-shell-capacity}
\end{equation}
Writing $x=i/T$, this becomes
\begin{equation}
 N_i\sim
 \frac{2\pi\lambda xT}{\sqrt{\lambda^2+p^2x^2}}.
 \label{eq:p-fold-shell-capacity-scaled}
\end{equation}
Therefore
\begin{equation}
 Q\sim q_p(\lambda)T^2,
 \label{eq:q-p-lambda-definition}
\end{equation}
where
\begin{align}
 q_p(\lambda)
 &=2\pi\lambda\int_0^1
 \frac{x\,dx}{\sqrt{\lambda^2+p^2x^2}} \\
 &=\frac{2\pi\lambda}
 {\lambda+\sqrt{\lambda^2+p^2}}.
 \label{eq:q-p-lambda}
\end{align}

\subsubsection{Length of one donut}

The major radius is
\begin{equation}
 R_0=2(1+\lambda)T.
 \label{eq:general-major-radius}
\end{equation}
A $p$-fold helix in shell $i$ has asymptotic length
\begin{equation}
 \ell_i
 \sim 4\pi T\sqrt{(1+\lambda)^2+p^2x^2}.
 \label{eq:p-fold-helix-length}
\end{equation}
Thus the length of one $T(pQ,Q)$ donut is
\begin{equation}
 L_p\sim A_p(\lambda)T^3,
 \label{eq:one-donut-asymptotic-length}
\end{equation}
with
\begin{equation}
 A_p(\lambda)
 =8\pi^2\lambda\int_0^1
 \frac{x\sqrt{(1+\lambda)^2+p^2x^2}}
 {\sqrt{\lambda^2+p^2x^2}}\,dx.
 \label{eq:A-p-lambda-integral}
\end{equation}

The integral has the closed form
\begin{equation}
 A_p(\lambda)
 =\frac{4\pi^2\lambda}{p^2}
 \left[
 \Phi_{p,\lambda}\left(\sqrt{\lambda^2+p^2}\right)
 -\Phi_{p,\lambda}(\lambda)
 \right],
 \label{eq:A-p-lambda-closed}
\end{equation}
where
\begin{equation}
 \Phi_{p,\lambda}(y)
 =y\sqrt{y^2+2\lambda+1}
 +(2\lambda+1)\operatorname{arsinh}
 \left(\frac{y}{\sqrt{2\lambda+1}}\right).
 \label{eq:Phi-p-lambda}
\end{equation}

For the double donut,
\begin{equation}
 L_{\mathrm{dbl}}\sim 2A_p(\lambda)T^3.
 \label{eq:double-donut-length}
\end{equation}
Since
\begin{equation}
 C_{\mathrm{dbl}}
 \sim 2(p+1)q_p(\lambda)^2T^4,
 \label{eq:double-donut-crossing-asymptotic}
\end{equation}
the asymptotic coefficient is
\begin{equation}
 \boxed{
 \alpha_p(\lambda)
 =\frac{2A_p(\lambda)}
 {\left[2(p+1)\right]^{3/4}q_p(\lambda)^{3/2}}.
 }
 \label{eq:alpha-p-lambda}
\end{equation}

For $p=1$ and $\lambda=1$, this gives
\begin{equation}
 \alpha_1(1)=11.640947\ldots,
 \label{eq:alpha-one-one}
\end{equation}
before the toroidal correction.  The construction remains conditional on
the no-overlap approximation.

\subsubsection{Numerical results}

The reported minimization over $\lambda\geq1$ gives
\begin{center}
\begin{tabular}{c|c|c|c}
 $p$ & $\lambda_{\mathrm{opt}}$ & $\alpha_p(1)$
     & optimized $\alpha_p$ \\
 \hline
 $1$ & $1$      & $11.640947$ & $11.640947$ \\
 $2$ & $1$      & $11.219985$ & $11.219985$ \\
 $3$ & $1.09418$      & $11.776016$ & $11.758138$ \\
 $4$ & $1.42033$      & $12.661548$ & $12.383066$
\end{tabular}
\end{center}

The optimizing values for $p=3,4$ are
\begin{equation}
 \lambda_{\mathrm{opt}}(3)\approx1.09418,
 \qquad
 \lambda_{\mathrm{opt}}(4)\approx1.42033.
 \label{eq:lambda-optima}
\end{equation}
Thus, under this double-donut shell model,
\begin{align}
 \alpha_{T(2Q,Q)\,\mathrm{double}}&\approx11.219985,\\
 \alpha_{T(3Q,Q)\,\mathrm{double}}&\approx11.758138,\\
 \alpha_{T(4Q,Q)\,\mathrm{double}}&\approx12.383066.
 \label{eq:double-donut-numerical-coefficients}
\end{align}

The $p=2$ construction improves the $p=1$ coefficient by approximately
\begin{equation}
 \frac{11.640947-11.219985}{11.640947}
 \approx3.62\%.
 \label{eq:p-two-improvement}
\end{equation}
The $p=3$ and $p=4$ constructions are less efficient than $p=2$, even after
optimizing the donut-hole radius.

These are asymptotic coefficients for the idealized straight-helix
construction.  They ignore toroidal lengthening, finite-shell effects, exact
integer shell populations, and the cone-piercing issue.  They therefore
represent conditional construction coefficients rather than rigorous
ropelength upper bounds.

\subsection{A more rigorous proof of the helical rectangle approximation}
\textbf{Author's note: In my paper, I made an approximation for the maximum number of helices that can fit in a cylinder, by finding the distance between parallel diagonals on a rectangle. I argued that the maximum number of helices permitted was slightly lower by a small $\epsilon<0.29$, and partially proved this by expanding the overlap constraint inequality to third order.}
\label{sec:exact-cylindrical-no-overlap}

The cubic argument can be replaced by an exact argument for the cylindrical
no-overlap constraint.  Define
\begin{equation}
 d_N^2(\theta)
 =2r^2\left[1-\cos\left(\theta+\frac{2\pi}{N}\right)\right]
 +h^2\theta^2,
 \qquad -\pi\leq\theta\leq\pi,
 \label{eq:exact-cylindrical-distance}
\end{equation}
where $r\geq2$ and $h\geq r$.  Let
\begin{equation}
 N_a=\frac{\pi hr}{\sqrt{h^2+r^2}},
 \qquad
 \varepsilon_0(r)
 =\pi\left(r-\frac{1}{\arcsin(1/r)}\right).
 \label{eq:Na-epsilon-zero}
\end{equation}

Then, for every real $N$ satisfying
\begin{equation}
 4\leq N\leq N_a-\varepsilon_0(r),
 \label{eq:exact-cylindrical-N-range}
\end{equation}
we have
\begin{equation}
 \boxed{
 d_N^2(\theta)\geq4
 \quad\text{for every }\theta\in[-\pi,\pi].
 }
 \label{eq:exact-cylindrical-separation}
\end{equation}
Consequently, every integer
\begin{equation}
 \boxed{
 4\leq N\leq
 \left\lfloor N_a-\varepsilon_0(r)\right\rfloor
 }
 \label{eq:integer-cylindrical-population}
\end{equation}
gives a rigorously non-overlapping cylindrical shell.

Equivalently, the correction parameter may be chosen in the interval
\begin{equation}
 \boxed{
 \varepsilon_0(r)\leq\varepsilon\leq N_a-4.
 }
 \label{eq:epsilon-fixed-parameter-range}
\end{equation}
The more conservative universal interval
\begin{equation}
 \varepsilon_0(r)\leq\varepsilon\leq2\pi-6
 \label{eq:epsilon-universal-range}
\end{equation}
is also valid in the stated regime.

\subsubsection{Proof}

Put
\begin{equation}
 n=N_a-\varepsilon,
 \qquad
 u=\frac{\pi}{n},
 \qquad
 \xi=\arcsin\frac1r.
 \label{eq:n-u-xi}
\end{equation}
Since $\varepsilon\geq\varepsilon_0(r)$ and $N_a\leq\pi r$,
\begin{equation}
 n\leq\frac{\pi}{\xi},
 \qquad	
 u\geq\xi.
 \label{eq:u-lower-bound}
\end{equation}
For the range $\varepsilon\leq2\pi-6$, together with $r\geq2$ and
$h\geq r$, one also has $n>4$, and therefore $u<\pi/4$.

At $\theta=0$,
\begin{equation}
 d_N^2(0)=4r^2\sin^2u\geq4,
 \label{eq:distance-at-zero}
\end{equation}
because $u\geq\xi$ and $r\sin\xi=1$.

Now consider negative $\theta$.  Write $\theta=-x$, with $x\geq0$.  For
$0\leq x\leq2u$,
\begin{equation}
 d_N^2(-x)
 =4r^2\sin^2\left(u-\frac{x}{2}\right)+h^2x^2.
 \label{eq:negative-theta-distance}
\end{equation}
Because $\sin$ is concave on $[0,\pi/4]$,
\begin{equation}
 \sin\left(u-\frac{x}{2}\right)
 \geq\frac{\sin u}{u}\left(u-\frac{x}{2}\right).
 \label{eq:concavity-sine-bound}
\end{equation}
Therefore
\begin{equation}
 d_N^2(-x)
 \geq
 \left(2r\sin u-\frac{r\sin u}{u}x\right)^2+h^2x^2.
 \label{eq:quadratic-lower-bound}
\end{equation}
The minimum of the quadratic expression on the right is
\begin{equation}
 \frac{4r^2h^2\sin^2u}
 {h^2+r^2(\sin u/u)^2}.
 \label{eq:quadratic-minimum}
\end{equation}
Thus it remains to prove
\begin{equation}
 r^2h^2\sin^2u
 \geq h^2+r^2\left(\frac{\sin u}{u}\right)^2.
 \label{eq:required-trigonometric-inequality}
\end{equation}

To establish this, define
\begin{equation}
 a(r)=r-\frac1\xi.
 \label{eq:a-of-r}
\end{equation}
The following trigonometric inequality holds for
$\xi\leq u\leq\pi/4$:
\begin{equation}
 \left(\frac1u+a(r)\right)
 \sqrt{1+r^2u^2-\left(\frac{u}{\sin u}\right)^2}
 \geq r.
 \label{eq:key-trigonometric-inequality}
\end{equation}

Indeed, let
\begin{equation}
 F(u)=\left(\frac1u+a(r)\right)
 \sqrt{1+r^2u^2-\left(\frac{u}{\sin u}\right)^2}.
 \label{eq:F-u}
\end{equation}
Differentiation shows that the sign of $F'(u)$ is the sign of
\begin{equation}
 S(u)=\frac{u^3\cos u}{\sin^3u}-1
 +a(r)u^3\left[r^2-
 \frac{\sin u-u\cos u}{\sin^3u}\right].
 \label{eq:S-u}
\end{equation}
For fixed $u$, the second term is increasing in $r$.  Its minimum occurs at
$r=\csc u$.  At that value,
\begin{equation}
 S(u)\geq
 \frac{u^4\cos u-\sin^4u}{\sin^4u}\geq0,
 \label{eq:S-nonnegative}
\end{equation}
where the final inequality follows from
\begin{equation}
 \frac{\sin u}{u}\leq(\cos u)^{1/4},
 \qquad 0\leq u\leq\frac\pi4.
 \label{eq:sine-cosine-inequality}
\end{equation}
Hence $F$ is increasing, and since $F(\xi)=r$, inequality
\eqref{eq:key-trigonometric-inequality} follows.

Because
\begin{equation}
 N_a=\frac\pi u+\varepsilon
 \geq\frac\pi u+\varepsilon_0(r),
 \label{eq:Na-lower-comparison}
\end{equation}
inequality \eqref{eq:key-trigonometric-inequality} implies
\begin{equation}
 N_a\sqrt{1+r^2u^2-\left(\frac{u}{\sin u}\right)^2}
 \geq\pi r.
 \label{eq:Na-square-root-bound}
\end{equation}
After substituting the definition of $N_a$ and simplifying, this is exactly
the required inequality
\eqref{eq:required-trigonometric-inequality}.  Therefore
\begin{equation}
 d_N^2(-x)\geq4,
 \qquad 0\leq x\leq2u.
 \label{eq:negative-region-separated}
\end{equation}

For $x\geq2u$,
\begin{equation}
 d_N^2(-x)\geq h^2x^2\geq h^2(2u)^2\geq4,
 \label{eq:negative-tail-separated}
\end{equation}
because $n\leq\pi r$, so $2u=2\pi/n\geq2/r$, and $h\geq r$.

The case $\theta\geq0$ is similar.  For
$0\leq\theta\leq\pi-2u$, both terms in $d_N^2(\theta)$ are no smaller
than their values at $\theta=0$, and hence
\begin{equation}
 d_N^2(\theta)\geq d_N^2(0)\geq4.
 \label{eq:positive-central-region}
\end{equation}
For $\theta\geq\pi-2u$, the $h^2\theta^2$ term alone exceeds $4$, since
$u<\pi/4$ and $h\geq2$.  Thus $d_N^2(\theta)\geq4$ throughout
$[-\pi,\pi]$, completing the proof.

\subsubsection{Numerical values of the correction}

The lower endpoint is
\begin{equation}
 \varepsilon_0(r)
 =\pi\left(r-\frac1{\arcsin(1/r)}\right).
 \label{eq:epsilon-zero-explicit}
\end{equation}
Representative values are
\begin{center}
\begin{tabular}{c|c}
 $r$ & $\varepsilon_0(r)$ \\
 \hline
 $2$   & $0.283185307$ \\
 $4$   & $0.133295257$ \\
 $6$   & $0.087963264$ \\
 $8$   & $0.065741949$ \\
 $10$  & $0.052508998$ \\
 $20$  & $0.026198507$ \\
 $100$ & $0.005236136$
\end{tabular}
\end{center}

The worst case is $r=2$, where
\begin{equation}
 \varepsilon_0(2)=2\pi-6=0.283185307.
 \label{eq:epsilon-worst-case}
\end{equation}
For large $r$,
\begin{equation}
 \varepsilon_0(r)=\frac\pi{6r}+O(r^{-3}).
 \label{eq:epsilon-asymptotic}
\end{equation}
Therefore the shell-count prescription can be made rigorous for the
cylindrical constraint by taking
\begin{equation}
 \boxed{
 N_{\max}=\left\lfloor
 \frac{\pi hr}{\sqrt{h^2+r^2}}
 -\pi\left(r-\frac1{\arcsin(1/r)}\right)
 \right\rfloor.
 }
 \label{eq:rigorous-cylindrical-shell-count}
\end{equation}

\subsubsection{Scope of the result}

This removes the cubic-approximation conjecture for the exact cylindrical
distance equation.  It proves that the shell occupancy used in the
construction is non-overlapping whenever $r\geq2$, $h\geq r$, and the
integer $N$ satisfies \eqref{eq:integer-cylindrical-population}.

It does not, by itself, make the complete toroidal upper bound rigorous.  A
separate toroidal correction is applied when the cylindrical construction is
bent into a torus; the non-overlap of that bent construction remains an
additional assumption.


\subsection{Fourth-order asymptotics of the toroidal correction}
\label{sec:toroidal-correction-fourth-order}

For the toroidal wrapping used in Klotz's correction, the excess over the
corresponding straight-helix length is asymptotically fourth order.  Let
$R$ be the major radius and $r$ the minor helix radius.  The relevant
toroidal and straight-helix speeds are proportional to
\begin{equation}
 v_T(\theta)=\sqrt{r^2+(R-r\cos\theta)^2},
 \qquad
 v_S=\sqrt{r^2+R^2}.
 \label{eq:toroidal-and-straight-speeds}
\end{equation}

Consequently, the correction factor is
\begin{equation}
 \mathcal C(\lambda)
 =\frac{1}{2\pi\sqrt{1+\lambda^{-2}}}
 \int_0^{2\pi}
 \sqrt{1-2\lambda^{-1}\cos\theta
 +\lambda^{-2}(1+\cos^2\theta)}\,d\theta,
 \label{eq:toroidal-correction-factor}
\end{equation}
where
\begin{equation}
 \lambda=\frac{R}{r}.
 \label{eq:major-minor-ratio}
\end{equation}

Set
\begin{equation}
 \varepsilon=\lambda^{-1}=\frac{r}{R}.
 \label{eq:epsilon-major-minor}
\end{equation}
For $\varepsilon<1$, a uniform expansion gives
\begin{equation}
 \frac{1}{2\pi}\int_0^{2\pi}
 \sqrt{1-2\varepsilon\cos\theta
 +\varepsilon^2(1+\cos^2\theta)}\,d\theta
 =1+\frac12\varepsilon^2
 +\frac18\varepsilon^4+O(\varepsilon^6).
 \label{eq:toroidal-average-expansion}
\end{equation}
The corresponding straight-helix factor is
\begin{equation}
 \sqrt{1+\varepsilon^2}
 =1+\frac12\varepsilon^2
 -\frac18\varepsilon^4+O(\varepsilon^6).
 \label{eq:straight-factor-expansion}
\end{equation}

Taking the ratio therefore yields
\begin{align}
 \mathcal C(\lambda)
 &=\frac{1+\frac12\varepsilon^2
          +\frac18\varepsilon^4+O(\varepsilon^6)}
         {1+\frac12\varepsilon^2
          -\frac18\varepsilon^4+O(\varepsilon^6)} \\[2pt]
 &=1+\frac14\varepsilon^4+O(\varepsilon^6).
 \label{eq:toroidal-ratio-expansion}
\end{align}
Thus
\begin{equation}
 \boxed{
 \mathcal C(\lambda)-1
 =\frac{1}{4\lambda^4}+O(\lambda^{-6})
 }
 \label{eq:toroidal-fourth-order}
\end{equation}
or, equivalently,
\begin{equation}
 \boxed{
 \lim_{\lambda\to\infty}
 \lambda^4\bigl(\mathcal C(\lambda)-1\bigr)=\frac14.
 }
 \label{eq:toroidal-fourth-order-limit}
\end{equation}

The cancellation of the $\lambda^{-2}$ term is the key point: both the
toroidal and straight-helical lengths acquire the same second-order
correction, leaving a fourth-order difference.

For example, when $R/r=2$, the asymptotic estimate gives
\begin{equation}
 1+\frac{1}{4(2)^4}=1.015625.
 \label{eq:toroidal-example-asymptotic}
\end{equation}
The exact integral gives approximately
\begin{equation}
 \mathcal C(2)=1.011229,
 \label{eq:toroidal-example-exact}
\end{equation}
which is consistent with the reported approximately $1.1\%$ correction for
the outermost helix.

One qualification is important.  In the full shell sum, the outer shells
retain $R/r\approx2$, so the total correction does not vanish as
$T\to\infty$; instead, it approaches the reported aggregate value
approximately $1.0042$.  The fourth-power law describes the correction for
a helix whose major-to-minor ratio itself tends to infinity.

\subsection{A polygonal improvement of the Wegner-based lower bound}
\label{sec:polygonal-wegner-improvement}

Let $N=Q-1$, and let $D$ be the convex hull of $N$ pairwise
non-overlapping unit disks.  Klotz's lower bound uses the Wegner area
estimate $A(D)\geq A_W(N)$ together with the ordinary isoperimetric
inequality,
\begin{equation}
 P(D)^2\geq 4\pi A(D),
 \label{eq:ordinary-isoperimetric-bound}
\end{equation}
which gives
\begin{equation}
 P(D)\geq W_N:=\sqrt{4\pi A_W(N)}.
 \label{eq:wegner-perimeter-bound}
\end{equation}
Cantarella's linking estimate then gives
\begin{equation}
 L_{T(Q,Q)}\geq Q\bigl(2\pi+P(D)\bigr).
 \label{eq:tqq-cantarella-bound}
\end{equation}

For the original form of Wegner's inequality, write
\begin{equation}
 A_W(N)=2\sqrt{3}\,(N-1)
 +(2-\sqrt{3})m_N+\pi,
 \qquad
 m_N=\left\lfloor\sqrt{12N-3}-3\right\rfloor.
 \label{eq:wegner-area}
\end{equation}

\subsubsection{A universal polygonal deficit}

Let $C$ denote the convex hull of the disk centers, and set
\begin{equation}
 p=P(C).
 \label{eq:center-hull-perimeter}
\end{equation}
Since $D=C+B_1$, Steiner's formulas give
\begin{equation}
 P(D)=p+2\pi,
 \qquad
 A(D)=A(C)+p+\pi.
 \label{eq:steiner-formulas}
\end{equation}

Suppose first that $C$ is a non-degenerate convex $m$-gon.  Consecutive
vertices of $C$ are distinct disk centers, so every side has length at least
$2$ and consequently
\begin{equation}
 p\geq2m.
 \label{eq:side-length-condition}
\end{equation}
The polygonal isoperimetric inequality gives
\begin{equation}
 A(C)\leq \frac{p^2}{4m}\cot\left(\frac{\pi}{m}\right).
 \label{eq:polygon-isoperimetric}
\end{equation}
Using $p\geq2m$ and
\begin{equation}
 \cot x\leq\frac1x-\frac{x}{3},
 \qquad 0<x\leq\frac{\pi}{3},
 \label{eq:cotangent-inequality}
\end{equation}
we obtain
\begin{equation}
 A(C)\leq \frac{p^2}{4\pi}-\frac{\pi}{3}.
 \label{eq:center-hull-deficit}
\end{equation}

For completeness, inequality \eqref{eq:cotangent-inequality} follows from
\[
 \tan x\geq\frac{3x}{3-x^2}.
\]
Indeed, if
\[
 f(x)=(3-x^2)\sin x-3x\cos x,
\]
then $f(0)=0$ and
\[
 f'(x)=x(\sin x-x\cos x)\geq0.
\]

If $C$ is a segment, then $A(C)=0$ and $p\geq4$, so
\eqref{eq:center-hull-deficit} remains valid.  Combining it with
\eqref{eq:steiner-formulas} yields the following universal inequality for
the convex hull of at least two non-overlapping unit disks:
\begin{equation}
 \boxed{
 A(D)\leq \frac{P(D)^2}{4\pi}-\frac{\pi}{3}.
 }
 \label{eq:universal-disk-hull-deficit}
\end{equation}

The term $\pi/3$ is the improvement over the circular isoperimetric bound.
It arises because the inner parallel body $C$ is a polygon whose sides are
bounded below by the disk diameter.

\subsubsection{Improved perimeter and ropelength bounds}

Combining \eqref{eq:universal-disk-hull-deficit} with
$A(D)\geq A_W(N)$ gives
\begin{equation}
 \boxed{
 P(D)\geq
 \sqrt{4\pi A_W(N)+\frac{4\pi^2}{3}}.
 }
 \label{eq:improved-wegner-perimeter}
\end{equation}
Therefore the improved lower bound for the $T(Q,Q)$ link is
\begin{equation}
 \boxed{
 L_{T(Q,Q)}
 \geq
 Q\left[
 2\pi+\sqrt{
 4\pi A_W(Q-1)+\frac{4\pi^2}{3}}
 \right].
 }
 \label{eq:improved-tqq-bound}
\end{equation}

Compared with Klotz's Wegner-based expression, the replacement is
\begin{equation}
 \sqrt{4\pi A_W(N)}
 \quad\longrightarrow\quad
 \sqrt{4\pi A_W(N)+\frac{4\pi^2}{3}}.
 \label{eq:wegner-replacement}
\end{equation}

A slightly sharper, implicit version can be obtained by retaining the exact
polygonal inequality.  Let $p_N$ be the least $p\geq4$ satisfying
\begin{equation}
 p+\pi+\frac{p}{2}\cot\left(\frac{2\pi}{p}\right)
 \geq A_W(N).
 \label{eq:implicit-polygonal-bound}
\end{equation}
Then
\begin{equation}
 P(D)\geq2\pi+p_N,
 \qquad
 L_{T(Q,Q)}\geq Q\bigl(4\pi+p_{Q-1}\bigr).
 \label{eq:implicit-tqq-bound}
\end{equation}
The explicit bound \eqref{eq:improved-wegner-perimeter} is slightly weaker but
has a closed form.

\subsubsection{Asymptotic size of the improvement}

Since
\begin{equation}
 A_W(N)=2\sqrt{3}\,N+O(\sqrt{N}),
 \label{eq:wegner-area-asymptotic}
\end{equation}
the two perimeter bounds have the same leading term.  Their difference is
\begin{align}
 \sqrt{4\pi A_W(N)+\frac{4\pi^2}{3}}
 -\sqrt{4\pi A_W(N)}
 &=\frac{4\pi^2/3}
 {\sqrt{4\pi A_W(N)+4\pi^2/3}+\sqrt{4\pi A_W(N)}} \\[2pt]
 &\sim
 \frac{2\pi^2}{3\sqrt{8\pi\sqrt{3}}}\,N^{-1/2}
 \approx0.99726\,N^{-1/2}.
 \label{eq:perimeter-improvement-asymptotic}
\end{align}

For $C=Q(Q-1)$, the corresponding improvement in the three-quarter
coefficient is therefore
\begin{equation}
 \alpha_{\mathrm{new}}(Q)-\alpha_{\mathrm{Wegner}}(Q)
 \sim\frac{0.99726}{Q}.
 \label{eq:alpha-improvement-asymptotic}
\end{equation}
The limiting coefficient itself is unchanged:
\begin{equation}
 \lim_{Q\to\infty}\alpha_{\mathrm{new}}(Q)
 =\sqrt{8\pi\sqrt{3}}
 \approx6.59782.
 \label{eq:new-limiting-coefficient}
\end{equation}

Representative values, using the original Wegner floor convention, are
\begin{center}
\begin{tabular}{c|c|c}
 $Q$ & $\alpha_{\mathrm{Wegner}}$ & $\alpha_{\mathrm{new}}$ \\
 \hline
 $7$    & $9.6887$ & $9.8554$ \\
 $10$   & $9.0909$ & $9.2010$ \\
 $20$   & $8.2883$ & $8.3402$ \\
 $100$  & $7.3240$ & $7.3340$ \\
 $1000$ & $6.8250$ & $6.8260$
\end{tabular}
\end{center}

For $Q=3,4,5,6$, the previously known exact small-packing bounds remain
stronger and should be combined with \eqref{eq:improved-tqq-bound} by taking
the maximum.



\section{A link version of the Rawdon--Simon Möbius--ropelength bound}
\label{sec:rawdon-link-bound}

This section records a link extension of the geometric argument of
Rawdon and Simon \cite{rawdon2002mobius}.  The principal point is that the proof does not require
the underlying one--manifold to be connected once the cross-energy is
included in the distal part of the integral.  The component count must,
however, be retained in the finite-length estimate.

\subsection{Conventions and definitions}

Let
\[
 \mathcal L=K_1\cup\cdots\cup K_m\subset\mathbb R^3
\]
be a smooth link with total length
\[
 \ell(\mathcal L)=\sum_{i=1}^m\ell(K_i)
\]
and thickness \(R(\mathcal L)\).  Its ropelength is
\[
 \operatorname{Rop}(\mathcal L)
 =
 \frac{\ell(\mathcal L)}{R(\mathcal L)}.
\]
After rescaling, we may assume \(R(\mathcal L)=1\), and we write
\[
 L=\ell(\mathcal L).
\]

For a component \(K_i\), let
\[
 E_{\mathrm s}(K_i)
 =
 \iint_{K_i\times K_i}
 \left(
 \frac{1}{|x-y|^2}
 -
 \frac{1}{\operatorname{arc}_{K_i}(x,y)^2}
 \right)\,ds_x\,ds_y
\]
denote the Rawdon--Simon/Freedman--He--Wang Möbius self-energy.  It is
normalized so that a round circle has energy \(4\).

There are two conventions for the cross-energy.  The once-counted
convention is
\[
 E_{\mathrm c}^{(1)}(K_i,K_j)
 =
 \iint_{K_i\times K_j}\frac{ds_x\,ds_y}{|x-y|^2},
 \qquad i<j,
\]
whereas the twice-counted convention is
\[
 E_{\mathrm c}^{(2)}(\mathcal L)
 =
 \sum_{i\ne j}
 \iint_{K_i\times K_j}
 \frac{ds_x\,ds_y}{|x-y|^2}
 =
 2\sum_{i<j}E_{\mathrm c}^{(1)}(K_i,K_j).
\]
Thus
\begin{align}
 E_{\mathrm M}^{(1)}(\mathcal L)
 &=\sum_i E_{\mathrm s}(K_i)
   +\sum_{i<j}E_{\mathrm c}^{(1)}(K_i,K_j), \label{eq:link-energy-once}\\
 E_{\mathrm M}^{(2)}(\mathcal L)
 &=\sum_i E_{\mathrm s}(K_i)
   +\sum_{i\ne j}E_{\mathrm c}^{(1)}(K_i,K_j). \label{eq:link-energy-twice}
\end{align}
For the Möbius-minimizing Hopf link, these conventions give respectively
\[
 \min E_{\mathrm c}^{(1)}=2\pi^2,
 \qquad
 \min E_{\mathrm c}^{(2)}=4\pi^2.
\]
The two conventions are therefore not interchangeable unless the
cross-term normalization is stated explicitly.

\subsection{The link estimate}

For each \(x\in K_i\), split the pairs involving \(x\) into:
\begin{enumerate}
 \item same-component proximal pairs with
 \(\operatorname{arc}_{K_i}(x,y)\leq\pi\);
 \item same-component distal pairs with
 \(\operatorname{arc}_{K_i}(x,y)\geq\pi\);
 \item all pairs \(y\in K_j\), \(j\ne i\).
\end{enumerate}
The first class is treated exactly as in the knot case.  By Schur's
comparison theorem and the circle calculation,
\begin{equation}
 E_{\mathrm{prox}}
 \leq \frac{2}{\pi}L.
\label{eq:link-proximal}
\end{equation}

Every component has length at least \(2\pi\), because the curvature is at
most one after the normalization \(R(\mathcal L)=1\), and Fenchel's
theorem gives total curvature at least \(2\pi\).  Hence for each
\(x\in K_i\), the set of same-component points at intrinsic distance
\(<\pi\) has length \(2\pi\).  The remaining same-component points have
length \(\ell(K_i)-2\pi\), while the other components have total length
\(L-\ell(K_i)\).  Consequently, the total distal length seen from \(x\)
is exactly
\[
 \bigl(\ell(K_i)-2\pi\bigr)+\bigl(L-\ell(K_i)\bigr)
 =L-2\pi.
\]

The thickness condition implies that every point in this distal set is at
Euclidean distance at least \(2\) from \(x\).  Moreover, the unit tube
around any distal subset lying in the ball \(B(x,r)\) lies in the shell
\[
 B(x,r+1)\setminus B(x,1).
\]
Since the tube is embedded and has volume \(\pi\) times its axial length,
the amount of distal curve that can lie within distance \(r\) is bounded
by
\begin{equation}
 h(r)=\frac{4}{3}\bigl((r+1)^3-1\bigr),
 \qquad r\geq2.
\label{eq:link-shell-capacity}
\end{equation}
In particular,
\[
 h(2)=\frac{104}{3}.
\]

Put \(M=L-2\pi\).  If \(M\leq104/3\), placing all distal length at the
smallest permitted radius gives
\[
 \Psi(M)=\frac{M}{4}.
\]
If \(M>104/3\), define
\begin{equation}
 P_M=\left(\frac34M+1\right)^{1/3}-1.
\label{eq:link-P}
\end{equation}
The shell-capacity argument and a Riemann-sum limit give
\begin{align}
 \Psi(M)
 &=
 \frac{104}{12}
 +\int_2^{P_M}\frac{4(r+1)^2}{r^2}\,dr \notag\\
 &=
 \frac83+4P_M-\frac{4}{P_M}
 +8\log\left(\frac{P_M}{2}\right).
\label{eq:link-psi}
\end{align}
Therefore the twice-counted link energy obeys
\begin{equation}
 \boxed{
 E_{\mathrm M}^{(2)}(\mathcal L)
 \leq
 L\,\Psi(L-2\pi)+4m .
 }
\label{eq:rawdon-link-explicit}
\end{equation}

The term \(4m\) comes from the regularizing arclength integrals on the
individual components:
\[
 \sum_{i=1}^m
 \ell(K_i)\left(
 2\int_{\pi}^{\ell(K_i)/2}\frac{dt}{t^2}
 \right)
 =
 \frac{2L}{\pi}-4m.
\]
Combining this identity with \eqref{eq:link-proximal} leaves \(4m\), just
as the knot proof leaves \(4\) when \(m=1\).

Since the once-counted energy has half the cross contribution,
\[
 E_{\mathrm M}^{(1)}(\mathcal L)
 =
 E_{\mathrm M}^{(2)}(\mathcal L)
 -
 \sum_{i<j}E_{\mathrm c}^{(1)}(K_i,K_j),
\]
we have the immediate, convention-independent consequence
\begin{equation}
 E_{\mathrm M}^{(1)}(\mathcal L)
 \leq E_{\mathrm M}^{(2)}(\mathcal L)
 \leq L\,\Psi(L-2\pi)+4m .
\label{eq:once-also-bounded}
\end{equation}
The sharper once-counted estimate would require repeating the shell
argument with only one copy of each cross pair; the displayed bound is the
safe common upper bound for both conventions.

\subsection{Power-law envelopes}

The explicit estimate \eqref{eq:rawdon-link-explicit} is preferable when
the number of components is known.  For a universal bound depending only
on \(L\), use \(m\leq L/(2\pi)\), again by Fenchel's theorem.  This gives
\begin{equation}
 E_{\mathrm M}^{(\nu)}(\mathcal L)
 \leq
 \Phi(L)
 :=
 L\,\Psi(L-2\pi)+\frac{2L}{\pi},
 \qquad \nu\in\{1,2\}.
\label{eq:universal-link-envelope}
\end{equation}
The leading term of \(\Psi\) is
\[
 \Psi(M)=4\left(\frac34M\right)^{1/3}
 +O(\log M),
\]
and hence
\begin{equation}
 \Phi(L)
 =
 4\left(\frac34\right)^{1/3}L^{4/3}
 +O(L\log L).
\label{eq:link-asymptotic}
\end{equation}
The asymptotic coefficient is
\[
 c_\infty
 =
 4\left(\frac34\right)^{1/3}
 =3.634241\ldots .
\]
Numerical one-variable maximization of the explicit function
\(\Phi(L)/L^{4/3}\) gives the conservative envelope
\begin{equation}
 \boxed{
 E_{\mathrm M}^{(\nu)}(\mathcal L)
 \leq c_{\mathrm{link}}L^{4/3},
 \qquad
 c_{\mathrm{link}}\simeq4.636 .
 }
\label{eq:link-c-link}
\end{equation}
The value \(4.636\) is the numerical envelope of the explicit link estimate
\eqref{eq:universal-link-envelope}; it is not the Rawdon--Simon knot
constant \(4.57\).  For a fully certified numerical statement, one should
round this value upward, for example to \(4.64\), after verifying the
one-variable maximum with interval arithmetic.  For a fixed number of
components, the sharper
\(m\)-dependent expression \eqref{eq:rawdon-link-explicit} has a smaller
envelope.  In particular, for \(m=1\), the original Rawdon--Simon
calculation recovers \(4.57\) (more precisely, approximately \(4.5626\))
and its asymptotic coefficient \(3.63\).

\subsection{Consequences for the Klotz energy--ropelength bound}

Klotz's Hopf-tambourine argument \cite{klotz2025geometric} uses \(N\) small circles linked to one
central circle.  Its crossing number is \(C=2N\), and in the
twice-counted \(4\pi^2\) convention the energy estimate is
\begin{equation}
 E_{\mathrm M}^{(2)}
 >
 4(N+1)+4\pi^2N
 =
 (2\pi^2+2)C+4.
\label{eq:tambourine-twice}
\end{equation}
Combining this with the universal link envelope
\eqref{eq:link-c-link} gives the rigorous conditional implication
\begin{equation}
 \operatorname{Rop}(\mathcal L)
 \gtrsim
 \left(
 \frac{2\pi^2+2}{c_{\mathrm{link}}}\,C
 \right)^{3/4}
 =
 3.19\,C^{3/4},
\label{eq:klotz-double-c-link}
\end{equation}
provided the Hopf tambourine is indeed the minimum-energy link at the
specified crossing number.  The \(+4\) in
\eqref{eq:tambourine-twice} can be retained for a finite-\(C\) version.

If one assumes that the sharper Rawdon--Simon numerical coefficient
\(4.57\) extends unchanged to the full twice-counted link energy, then
\eqref{eq:klotz-double-c-link} becomes
\[
 \operatorname{Rop}(\mathcal L)
 \gtrsim
 \left(
 \frac{2\pi^2+2}{4.57}\,C
 \right)^{3/4}
 =
 3.22\,C^{3/4}.
\]
This is Klotz's quoted coefficient.  The explicit link proof above
justifies the same mechanism but, without an additional optimization of
the component-count term, gives the slightly more conservative \(3.19\).
Using only the asymptotic coefficient \(c_\infty=3.634241\ldots\) gives
\[
 \operatorname{Rop}(\mathcal L)
 \gtrsim 3.82\,C^{3/4}
\]
at the level of leading asymptotics.

\subsection{The factor-of-two convention}

\textbf{Author's note: there is a choice of convention when computing the Mobius energy of links regarding whether each pair of line elements contribute once or twice to the energy, such that the minimum Mobius energy of the Hopf link is either $2\pi^2$ or $4\pi^2$. The AI was asked which was more likely to be correct.}

Under the once-counted convention, the same tambourine has
\begin{equation}
 E_{\mathrm M}^{(1)}
 >
 4(N+1)+2\pi^2N
 =
 (\pi^2+2)C+4.
\label{eq:tambourine-once}
\end{equation}
Using the same universal link envelope gives
\begin{equation}
 \operatorname{Rop}(\mathcal L)
 \gtrsim
 \left(
 \frac{\pi^2+2}{c_{\mathrm{link}}}\,C
 \right)^{3/4}
 =
 2.02\,C^{3/4}.
\label{eq:klotz-once-c-link}
\end{equation}
At the level of the large-\(L\) asymptotic coefficient
\(c_\infty=3.634241\ldots\), the corresponding once-counted coefficient is
\[
 \left(\frac{\pi^2+2}{c_\infty}\right)^{3/4}
 =2.43\ldots .
\]
The factor of two belongs to the definition of the cross term, not to the
ropelength inequality itself.  The self-energy is unchanged, so the full
link energies are not related by a global factor of two:
\[
 E_{\mathrm M}^{(2)}
 =
 E_{\mathrm M}^{(1)}
 +\sum_{i<j}E_{\mathrm c}^{(1)}(K_i,K_j).
\]
For the tambourine, this changes the leading energy-per-crossing
coefficient from \(\pi^2+2\) to \(2\pi^2+2\), and therefore changes the
ropelength coefficient by the \(3/4\) power of the corresponding ratio.
Consequently, the numerical lower bound must always state both the
cross-energy convention and the convention used in the Rawdon--Simon
constant.

\subsubsection*{Justification of Proper Convention}

Let \(K_1,\ldots,K_m\) be the components of a link. Denote the
self-energy of \(K_i\) by
\[
E_i=E(K_i),
\]
and define the Möbius cross-energy between distinct components by
\[
X_{ij}
=
\iint_{K_i\times K_j}
\frac{ds_i\,ds_j}{\lvert x-y\rvert^2}.
\]

Because
\[
X_{ij}=X_{ji},
\]
the interaction between two components should be counted only once.
The total link energy used by Freedman, He, and Wang is therefore
\[
\boxed{
E_{\mathrm{link}}(K_1,\ldots,K_m)
=
\sum_{i=1}^{m}E_i
+
\sum_{1\leq i<j\leq m}X_{ij}
}
\]
or, equivalently,
\[
\boxed{
E_{\mathrm{link}}
=
\sum_{i=1}^{m}E_i
+
\frac12
\sum_{\substack{i,j=1\\i\neq j}}^{m}X_{ij}.
}
\]

The factor \(1/2\) removes the double counting in the ordered sum.
This is the convention introduced for link energies by
Freedman, He, and Wang~\cite{mobius}.

An alternative convention is
\[
E_{\mathrm{ordered}}
=
\sum_{i=1}^{m}E_i
+
\sum_{\substack{i,j=1\\i\neq j}}^{m}X_{ij}.
\]
For a two-component link, this gives
\[
E_{\mathrm{ordered}}=E_1+E_2+2X_{12},
\]
whereas the single-count convention gives
\[
E_{\mathrm{link}}=E_1+E_2+X_{12}.
\]

The single-count convention is the more natural total interaction
energy for three reasons. First, each physical pair of points is counted
once. Second, it agrees with the usual definition of Möbius cross-energy
for a two-component link. Third, when components are separated by a
large distance, \(X_{ij}\to0\), so the energy approaches the sum of the
individual self-energies.

For the standard Hopf-link configuration, the cross-energy in this
normalization is
\[
X_{12}=2\pi^2.
\]
Using the Freedman--He--Wang normalization, a round circle has
self-energy \(4\), so the total energy of the Hopf link is
\[
E_{\mathrm{Hopf}}=4+4+2\pi^2
=
8+2\pi^2.
\]
The ordered convention would instead assign the value
\[
8+4\pi^2,
\]
because it counts the cross-interaction twice.

Consequently, the natural link extension of a self-energy--ropelength
inequality such as the Rawdon--Simon inequality is
\[
E_{\mathrm{link}}
=
\sum_i E(K_i)
+
\sum_{i<j}E_{\times}(K_i,K_j).
\]
If the ordered convention is used instead, every cross-energy coefficient
in the corresponding ropelength inequality must also be doubled.
The two conventions are mathematically equivalent only when all
minimum-energy values and inequality constants are converted
consistently.

\subsection{What is proved and what remains conditional}

The geometric estimate \eqref{eq:rawdon-link-explicit} is the rigorous
link modification: it follows from curvature, thickness, tube-volume
packing, and the same shell rearrangement used by Rawdon and Simon.  The
universal envelope \eqref{eq:link-c-link} is a convenient consequence.

The following separate assertion is not proved here:
\[
 \text{the Hopf tambourine minimizes Möbius energy among all links
 with a prescribed crossing number.}
\]
Thus the coefficients \(3.19\), \(3.22\), \(3.82\), \(2.02\), and \(2.43\) are
respectively conditional crossing-number lower-bound coefficients: they
become universal only after an independent lower bound for the minimum
Möbius energy at crossing number \(C\) has been established.

\textbf{Author's note: during this query path, ChatGPT got confused between two different inequalities in the Rawdon-Simon paper and falsely reported a ropelength bound of $2.135C^{3/4}$ and had to be corrected. This should serve as another warning for anyone taking LLM outputs purely at face value.}

\end{document}